\documentclass{article}
\usepackage{a4wide}

\usepackage{graphicx}
\usepackage{amssymb}
\usepackage{amsmath}
\usepackage{graphicx} 
\usepackage{geometry}
\usepackage[dvipsnames]{xcolor}
\usepackage{mathrsfs} 
\usepackage{bm}             
\usepackage{enumitem}
\setlist[enumerate]{leftmargin=*,noitemsep, topsep=3pt,parsep=0pt,partopsep=0pt}
\setlist[itemize]{leftmargin=*,noitemsep, topsep=3pt,parsep=0pt,partopsep=0pt}
\usepackage{setspace}
\usepackage{rotating}

\usepackage{numprint}
\usepackage{algorithmic}
\usepackage{subcaption}
\usepackage[linesnumbered,ruled]{algorithm2e}

\usepackage[scientific-notation=true]{siunitx} 
\newcommand{\numf}[1]{\num[round-mode=places, scientific-notation=false, round-precision=0]{#1}}

\usepackage{jabbrv}

\newcommand{\smaz}[1]{{\color{brown}{ }}}

\usepackage{hyperref}
\hypersetup{
colorlinks   = true, 
  urlcolor     = blue, 
  linkcolor    = blue, 
  citecolor   = red, 
  pdftitle={Hybrid Schwarz preconditioners for linear systems arising from
  $hp$-discontinuous Galerkin method},
  pdfauthor={V.  Dolej{\v s}{\'\i}, and T. Hammerbauer}
}

\newcommand{\vertical}[1]{\begin{sideways}{#1}\end{sideways}}

\renewcommand\em {\it}

\DeclareMathOperator*{\supp}{{supp}}
\DeclareMathOperator*{\diag}{{diag}}

\newcommand{\range}{\mathrm{Range\,}}

\newtheorem{theorem}{Theorem}[section]
\newtheorem{lemma}[theorem]{Lemma}
\newtheorem{definition}[theorem]{Definition}
\newtheorem{remark}[theorem]{Remark}

\newtheorem{assumption}[theorem]{Assumption}
 \newtheorem{corollary}[theorem]{Corollary}

\usepackage{thmtools}   
\newcommand\bifont {\bm}

\newcommand\krt{{\mathscr T}}

\newcommand\Th{{\krt_h}}
\newcommand\Thi{{\krt_{h,i}}}
\newcommand\ThH{{\krt_H}}
\newcommand\ThHs{{\krt_H^{\#}}}

\newcommand\Fh{{\Gamma_h}}
\newcommand\FhI{{\Gamma_h^I}}
\newcommand\FhD{{\Gamma_h^B}}

\newcommand\FhB{{\Gamma_h^B}}
\newcommand\FhID{{\Gamma_h}}

\newcommand\diam{{\rm diam}}

\newcommand\bke{\mbox{\boldmath$e$\unboldmath}}

\newcommand\bkg{\mbox{\boldmath$g$\unboldmath}}

\newcommand\bkn{\mbox{\boldmath$n$\unboldmath}}

\newcommand\bku{\mbox{\boldmath$u$\unboldmath}}
\newcommand\bkv{\mbox{\boldmath$v$\unboldmath}}
\newcommand\bkw{\mbox{\boldmath$w$\unboldmath}}
\newcommand\bkx{\mbox{\boldmath$x$\unboldmath}}
\newcommand\bky{\mbox{\boldmath$y$\unboldmath}}
\newcommand\bkz{\mbox{\boldmath$z$\unboldmath}}

\newcommand\bfzero {\bifont{0}}
\newcommand\mA{\bifont{A}}

\newcommand\mR{\bifont{R}}
\newcommand\mP{\bifont{P}}
\newcommand\mI{\bifont{I}}
\newcommand\mN{\bifont{N}}

\newcommand\Ri{R_i}
\newcommand\RiT{R_i^\T}
\newcommand\RjT{R_j^\T}
\newcommand\mRi{\mR_i}
\newcommand\mRiT{\mR_i^\T}
\newcommand\mRjT{\mR_j^\T}
\newcommand\tPi{\tilde{P}_i}
\newcommand\mtPi{\tilde{\mP}_i}
\newcommand\PPi{P_i}
\newcommand\PPj{P_j}
\newcommand\PPh{\hat{P}}
\newcommand\mPP{\mP}
\newcommand\mPPi{\mPP_i}

\newcommand\RN{R_0}
\newcommand\RNT{R_0^\T}
\newcommand\mRN{\mR_0}
\newcommand\mRNT{\mR_0^\T}

\newcommand\PPN{P_0}
\newcommand\mPPN{\mP_0}

\newcommand\add{\mathrm{add},2}

\newcommand\Padd{P_{\add}}

\newcommand\mPadd{\mP_{\add}}
\newcommand\mPad{\hat{\mP}}

\newcommand\mNadd{\mN_{\add}}
\newcommand\mNad{\hat{\mN}}
\newcommand\mNN{\mN_0}

\newcommand\hyb{\mathrm{hyb}}

\newcommand\mNhyb{\mN_{\hyb}}

\newcommand\mPhyb{\mP_{\hyb}}
\newcommand\Phyb{P_{\hyb}}

\newcommand\R{{\mathbb R}}

\newcommand\hK{{h_K}}
\newcommand\hKp{{h_{K'}}}
\newcommand\nng{\bkn_{\gamma}}
\newcommand\hg{h_\gamma}
\newcommand\pg{p_\gamma}

\newcommand\pd {{\partial}}
\newcommand\pdK {\partial K}

\newcommand\dx{\,{\rm d}x}

\newcommand\dS{\,{\rm d}S}

\newcommand\Om{\Omega}
\newcommand\oOm  {\overline{\Om}}
\newcommand\gom{{\Gamma}}
\newcommand\gomN{{\Gamma_N}}
\newcommand\gomD{{\Gamma_D}}

\newcommand\gomij{{\Gamma_{ij}}}

\newcommand\Ldo{L^2(\Om)}
\newcommand\Ldoi{L^2(\Om_i)}

\newcommand\Ldd{L^2(\gamma)}

\newcommand\norm[2]{{\left\|{#1}\right\|_{#2}^{}} }

\newcommand\normP[3]{{\left\|{#1}\right\|_{#2}^{#3}} }

\newcommand\bnormP[3]{{\big\|{#1}\big\|_{#2}^{#3}} }

\newcommand\llbracket {[\![}
\newcommand\rrbracket {]\!]}

\newcommand\jump[1]{\llbracket{#1}\rrbracket}

\newcommand\aver[1]{\left\langle{#1}\right\rangle}

\newcommand\Lsp[3]{{( {#1}, {#2} )_{#3} }}
\newcommand\LspK[2]{{( {#1}, {#2} )_{K} }}
\newcommand\Lspg[2]{{( {#1}, {#2} )_{\gamma} }}

\newcommand\nolim {}

\newcommand\sumK  {{ \sum\nolim_{K\in\Th}     }}
\newcommand\sumF  {{ \sum\nolim_{\gamma\in\Fh} }}

\newcommand\sumFij {{ \sum\nolim_{\gamma\in\gomij} }}
\newcommand\sumFID{{ \sum\nolim_{\gamma\in\FhID} }}

\newcommand\sumFD {{ \sum\nolim_{\gamma\in\FhD} }}

\newcommand\qed{{\hfill$\Box$} }

\newcommand\dof {n}

\newcommand\vi {v_i}
\newcommand\vj {v_j}
\newcommand\tvi {\tilde{v}_i}
\newcommand\tvj {\tilde{v}_j}

\newcommand\K {a} 
\newcommand\barK {\overline{\K}} 
\newcommand\ubarK {\underline{\K}} 

\newcommand\cK {\mathcal{K}}

\newcommand\cA {\mathcal{A}}
\newcommand\Ah {\cA_h}

\newcommand\Ahi {\cA_{h,i}}

\newcommand\Bh {\mathcal{B}_h}

\newcommand\lh {{g_h}}

\newcommand\T{{\intercal}}

\newcommand{\DoF} {{\mathrm{DoF}}}

\newcommand\ve {\varepsilon}
\newcommand{\bve}{\bm\varepsilon}

\newcommand{\vp} {{\varphi}}

\newcommand\Shp {S_{h}}
\newcommand\Ship {S_{h,i}}
\newcommand\Shjp {S_{h,j}}
\newcommand\SHP {S_{h,0}}
\newcommand\SHHP {S_{H}}

\newcommand\TOL {{{\omega}}}

\def \u {{u}} 
 
\def \uh {u_h}

\def \wh {w_h}
\def \wH {w_H}

\newcommand\restr[2]{{
  \left.\kern-\nulldelimiterspace 
  #1 
  \vphantom{\big|} 
  \right|_{#2} 
  }}

\newcommand\pK{p_K}
\newcommand\pKp{p_{K'}}
\newcommand\qK{q_{\cK}}
\newcommand\dK{{\pd K}}
\newcommand\dcK{{\pd \cK}}

\newcommand\NS {N_S}
\newcommand\rrel {r_{\mathrm{rel}}}

\newcommand\fac{\mathrm{fac}}
\newcommand\ass{\mathrm{ass}}
\newcommand\fl {\mathsf{fl}}
\newcommand\flfac {\fl_{\fac}}
\newcommand\flass {\fl_{\ass}}
\newcommand\FFfac {\mathsf{Fl}_{\fac}}
\newcommand\FFass {\mathsf{Fl}_{\ass}}
\newcommand\Fl {\mathsf{flops}}
\newcommand\MFl {\mathsf{Mflops}}

\newcommand\iter {\mathsf{iter}}
\newcommand\iterN {\iter_{\mathrm{N}}}
\newcommand\iterL {\iter_{\mathrm{L}}}
\newcommand\comm {\mathsf{com}}
\newcommand\Mcomm {\mathsf{Mcom}}

\newcommand{\LdcK}{L^2(\dcK)}
\newcommand{\LdK}{L^2(\dK)}
\newcommand{\dg}{\vert \! \vert \! \vert}
\newcommand{\dgI}{\vert \! \vert \! \vert_{\sigma}}

\newcommand{\CW}{C_W}
\newcommand{\Coo}{C_0}
\newcommand{\Co}{C_0^2}
\newcommand{\vh}{v_h}

\newcommand{\sumfh}{\sumF}

\newcommand{\sbar}{\overline{\sigma}}

\begin{document}

\title{Hybrid Schwarz preconditioners for linear systems arising from
  $hp$-discontinuous Galerkin method\thanks{This work has been supported by the
    Czech Science Foundation grant no. 20-01074S (V.D.),
the Charles University grant SVV-2025-260711 (T.H.).
V.D. acknowledges the membership in the Ne{\v c}as Center for Mathematical Modeling ncmm.karlin.mff.cuni.cz.
  }
}


\author{V{\'\i}t Dolej{\v s}{\'\i}, Tom{\'a}{\v s} Hammerbauer\footnote{V. Dolej{\v s}{\'\i}, T. Hammerbauer,
  Charles University, Faculty of Mathematics and Physics, Sokolovsk{\'a} 83, Prague,
Czech Republic,
              {vit.dolejsi@matfyz.cuni.cz, hammerbt@karlin.mff.cuni.cz}            
}
}




\maketitle

\begin{abstract}
  We deal with the numerical solution of linear elliptic problems with varying diffusion
  coefficient by the $hp$-discontinuous Galerkin method.
  We develop a two-level hybrid Schwarz preconditioner for the arising linear algebraic systems.
  The preconditioner is additive with respect to the local components and multiplicative
  with respect to the mesh levels. 
  We derive the $hp$ spectral bound of the preconditioned operator in the form $O((H/h)(p^2/q))$,
  where $H$ and $h$ are the element sizes
  of the coarse and fine meshes, respectively, and $p$ and $q$ are the polynomial approximation
  degrees on the fine and coarse meshes. Further, we present a numerical study
  comparing the hybrid Schwarz preconditioner with the standard additive one
  from the point of view of the speed of convergence and also computational costs.
  Moreover, we investigate the convergence of both techniques
  with respect to the diffusivity variation and
  to the domain decomposition (non-)respecting the material interfaces.
  Finally, the combination with
  a $hp$-mesh adaptation for the solution of nonlinear problem demonstrates the potential
  of this approach.
\end{abstract}



\section{Introduction}
\label{sec:intro}

The discontinuous Galerkin (DG) method \cite{PietroErn2012,DGM-book}
exhibits a powerful technique for the
numerical solution of partial differential equations. 
In particular, the piecewise polynomial discontinuous approximation is well suited
for $hp$-mesh adaptation which exhibits an excellent strategy for problems with
local singularities, material interfaces, and it
gives, under some assumption, an exponential rate of convergence.
On the other hand, the DG discretization leads to larger (but sparser)  algebraic systems
in comparison to conforming finite element methods.

Therefore, a development of efficient iterative solvers is demanding.
Among many techniques,
a prominent role is played by the {\em domain decomposition} (DD) methods
\cite{QuarteroniValli00,Nataf-DDM,ToselliWidlund-DD05,Gander_ETNA08} that can use the power of multiprocessor computers.
The main idea is to divide the global problem into several subproblems
(based on the decomposition of the computational domain) and solve them
independently. The transfer of information among the subproblems
can be accelerated by considering a global (coarse) problem.
The more frequent approach is the use of DD techniques as preconditioners
for Krylov space iterative methods, such as the conjugate gradient method.

The domain decomposition techniques
for discontinuous Galerkin approximations were considered
in many works, let us mention \cite{Karashian2001,AntoniettiAyuso_M2NA07,AntoniettiALL_JSC14,PaznerKolev_CAMC22,GopalakrishnanKanschat_23},
a further paper \cite{AntoniettiAyuso_CCP09} related to the super-penalty DG technique, 
works \cite{GanderHajian_SINUM14,GanderHajian_MC18,BarrenecheaALL_CMAM19}
dealing with hybridizable DG methods, and
papers
\cite{Canuto-2013-BPC,Kim-2014-BAC,Tu-2021-BAA,Dryja-2007-DGD}
related to
balancing domain decomposition by constraints (BDDC) variants of domain decomposition.
In particular, we mention papers
\cite{AntoniettiHouston_JSC11,AntoniettiALL_IJNAM16} that analyze two-level
nonoverlapping additive and multiplicative Schwarz preconditioners.
See also \cite{DryjaKrzyzanowski_NM16,Krzyzanowski_NMPDE16,AntoniettiALL_MC20}
containing an alternative approach when the
coarse mesh is created by several subdomains treated in parallel. 

In contrary to conforming finite element methods (FEM), the advantage of
nonoverlapping discontinuous Galerkin domain decomposition approach is
that each degree of freedom belongs to only one subdomain,
any specific operator at subdomain interfaces need not be
constructed, and the coarse operator is defined using the same variational formulation
as the original problem.

In this paper, we introduce the Schwarz {\em two-level nonoverlapping hybrid Schwarz} technique
proposed in \cite{Mandel_94} for FEM,
when the preconditioner operator is additive with respect to the local components and multiplicative
with respect to the mesh levels. 
The hybrid Schwarz technique was used in various applications, for example
\cite{BarkerCai_SISC10,Scacchi_CMAME08,HeinleinLanser_SISC20}.

We present the discretization of a linear elliptic problem using the symmetric
interior penalty Galerkin (SIPG) method. Assuming a suitable indexing of mesh elements,
the one-level additive Schwarz method is equivalent to the block Jacobi iterative method.
Then we introduce the two-level additive and hybrid Schwarz methods, and adapting the techniques
from \cite{AntoniettiALL_IJNAM16}, we derive the $hp$ spectral bound of the
two-level nonoverlapping hybrid Schwarz preconditioner, which is the first novelty of this paper.
We note that a minor modification of the abstract theory from
\cite[Theorem~2.13]{ToselliWidlund-DD05} is required, cf. Remark~\ref{rem:Widlund}.

Moreover, the core of this paper is the detailed numerical study of the convergence of the conjugate
gradient (CG) method with the two-level additive and hybrid Schwarz preconditioners.
The aim is to achieve the {\em weak scalability} of the solver, namely
the number of iterations of the algebraic solver does not increase when the size
of sub-problems is kept fixed. Additionally, we compare both preconditioners from the
point of view of the number of floating point operations and number of communications among
the computer cores. Furthermore, we consider a problem with piecewise constant diffusion
and investigate the dependence of the convergence on the diffusivity variation and
the domain decomposition construction, which does or does not respect the material interfaces.

Finally, the potential of both techniques is demonstrated by numerical solution of a nonlinear elliptic equation
in combination with a mesh adaptation. In particular, the arising nonlinear algebraic system is
solved by the Newton method and at each Newton step, the corresponding linear system is solved
by CG with the presented preconditioners. The mesh adaptation is carried out by
the {\em anisotropic $hp$-mesh adaptation} technique from \cite{AMA-book}, which serves as a test
of the robustness of the hybrid and additive preconditioners with respect to $hp$-adaptation
and the anisotropy of the mesh elements.

The rest of the paper is the following.
In Section~\ref{sec:DGM}, we introduce the SIPG discretization of the model problem.
In Section~\ref{sec:DDM}, we introduce the
two-level hybrid Schwarz preconditioner and in Section~\ref{sec:anal} we derive the corresponding
spectral bound. The numerical study demonstrating the performance of additive and hybrid preconditioners is given in Section~\ref{sec:numer}.
Several concluding remarks are given in Section~\ref{sec:concl}.

\section{Problem definition}
\label{sec:DGM}

We use the standard notation for the Lebesgue space $L^2(\Om)$, the Sobolev space $H^{k}(\Om) = W^{k,2}(\Om)$, and $P^p(M)$ is the space of polynomials of degree at most $p$ on
$M \subset \R^d$.  
Moreover, we use the notation
\begin{align}
  \label{Lsp}
  \Lsp{u}{v}{M}&=\int_M u\, v \dx\qquad \mbox{ for } M \subset \R^s,\ s=1,\dots, d.
\end{align}

Let $\Om \in \R^d$, $d=2,3$ be a bounded convex
polygonal domain with Lipschitz boundary $\gom:=\pd\Om$.
Let $f\in L^2(\Om)$ be a source term and
$\K\in L^\infty(\Om)$ be a diffusion such that
\begin{align}\label{Kij}
  0 < \ubarK \leq \K(x) \leq \barK \quad \forall x\in\Om,
  \qquad \mbox{where} \quad \ubarK,\barK\in\R.
\end{align}
Later, we introduce an additional requirement for the diffusion function $\K$,
cf.~Assumption~\ref{ass:K}.
We consider the linear model problem with the Dirichlet boundary condition
\begin{subequations}
\label{prob1}
  \begin{align}
  \label{prob1a}
  - \nabla\cdot(\K \nabla u) 
  &= f \ \ \qquad \text{in } \Om, \\ 
  \label{prob1b}
  u &= u_D  \qquad\text{on } \gom,
  \end{align}
 \end{subequations}
where $u: \Om \to \R$ is an unknown scalar function defined in $\Om$ and
$u_D\in L^2(\gomD)$. 
Obviously, problem \eqref{prob1} has a unique weak solution.
We note that the assumption of convexity of $\Om$ can be relaxed by the technique
from \cite{AntoniettiALL_MC20}.

\subsection{Discontinuous Galerkin discretization} 

Let $\Th$ be a partition of $\oOm$ consisting of a finite
number of $d$-dimensional simplices $K$ with mutually
disjoint interiors. The symbol $\dK$ denotes a boundary of $K \in \Th$,
and $h_K = \diam(K)$ is its diameter.
The approximate solution of \eqref{prob1} is sought in
the finite-dimensional space 
\begin{align}
  \label{Shp}
  \Shp  = \left\{ v \in \Ldo; \, \restr{v}{K} \in P^{\pK}(K) \, \forall K \in \Th \right\}, 
\end{align}
where $\pK$ is the local polynomial degree assigned to each $K\in\Th$. 
We assume that the ratio of polynomial degrees $\pK$ and $\pKp$  of any pair of elements
$K$ and $K'$
that share a face (edge for $d=2$) is bounded.

\begin{assumption}
  \label{ass:K}
  We assume that the diffusion function $\K$ in \eqref{Kij} has
  traces on $\dK$, $K\in\Th$ such that
  $\K\in L^\infty(\dK)$, $K\in\Th$.
  This is guaranteed when $\K$ is piecewise regular.
\end{assumption}

By $\Fh$ we denote the union of all faces $\gamma$ contained in $\Th$.
Moreover,
$\FhI$ and $\FhB$ denote the union of interior and boundary faces of $\Fh$, respectively.
For each $\gamma\in\Fh$, we consider a unit normal vector $\nng$,
its orientation can be chosen arbitrarily for the interior faces.
Symbols $\jump{v}_\gamma$ and $\aver{v}_\gamma$ denote
the jump of $v$ multiplied by $\nng$  and the mean value of $v\in\Shp$
on $\gamma\in\FhI$, respectively.
For $\gamma\in\FhB$, we put  $\jump{v}_\gamma = v\bkn_\gamma$ and $\aver{v}_\gamma=v$.
If there is no risk of misunderstanding
we drop the subscripts $_\gamma$.


Let $\gamma\in \FhI$ and
$\gamma \subset \pdK\cap \pdK'$, $K\not=K'$, we set
\begin{align}
  \label{gamma}
  \hg = \max(\hK, \hKp), \quad 
  \pg = \max(\pK, \pKp),\quad
  \K_\gamma(x)= \max(\K(x)|_{\pdK}, \K(x)|_{\pdK'}),\ x\in\gamma,
\end{align}
where $\K(x)|_{\pdK}$ and $\K(x)|_{\pdK'}$ denote  the traces of $\K$
  on $\gamma\subset \pdK\cap \pdK'$, which exist due to  Assumption~\ref{ass:K}.
For $\gamma\subset \pdK \cap\FhB$, we set
\begin{align}
  \label{gammaB}
  \hg=\hK,\quad  \pg = \pK,\quad 
\K_\gamma(x)= \K(x)|_{\pdK},\ x\in\gamma.
\end{align}

Using \eqref{Lsp}, we define the forms for $u,v \in \Shp$,
(cf. \cite[Section~4.6]{DGM-book} or \cite{Houston-book})
\begin{align} 
  \Ah(u,v):= &  \sumK\LspK{\K \nabla u}{ \nabla v } 
   -  \sumFID\Big( \Lspg{ \aver{\K\nabla \u} -\sigma\jump{\u}}{\jump{v}}
  +  \Lspg{ \aver{\K\nabla v} }{\jump{u}} 
  \Big), \notag  \\  
  \label{Ah}
  \lh( v ) := & \Lsp{f}{v}{\Om}
  - \sumFD \Big(\Lspg{ \K \nabla v }{u_D\,\nng} - \Lspg{ \sigma u_D} { v }\Big),
\end{align} 
where
the penalty parameter $\sigma>0$
is chosen by
\begin{align}
  \label{sigma}
  \restr{\sigma}{\gamma} 
  = {\CW \K_\gamma\, \pg^2}/{\hg},\ \gamma \in \FhID,
\end{align}
with
a constant $\CW > 0$, chosen sufficiently large to guarantee convergence of the method,
see \cite[Chapter~2]{DGM-book}.  Then the discrete problem reads.
\begin{definition}
  We say that $\uh \in \Shp$ is the {\em approximate solution} of \eqref{prob1}
  by {\em symmetric interior penalty Galerkin (SIPG)} method
  if 
\begin{align} \label{DG}
 \Ah(\uh,v_h) = \lh(v_h)\quad \forall v_h \in \Shp.
\end{align} 
\end{definition}

\section{Domain decomposition method}
\label{sec:DDM}

We follow the approach of
\cite{AntoniettiHouston_JSC11,AntoniettiALL_IJNAM16}, where more details can be found.
We consider a nonoverlapping domain decomposition of $\Om$ as a set of open subdomains
$\Om_i,\ i=1,\dots,N$ such that $\Om_i\cap\Om_j=\emptyset$ and $\oOm = \cup_{i=1,\dots, N}\oOm_i$.
We denote by $\gom_i$ the boundaries of $\Om_i$, $i=1,\dots,N$.
The subdomains $\Om_i$ are constructed as a union of some elements $K\in\Th$,
the corresponding meshes are denoted as $\Thi,\ i=1,\dots, N$.
We consider finite-dimensional spaces
\begin{align}
  \label{Ship}
  \Ship  = \left\{ v \in \Ldoi; \, \restr{v}{K} \in P^{\pK}(K) \, \forall K \in \Thi \right\},
  \qquad i=1,\dots, N.
\end{align}

In addition to the mesh $\Th$, we consider a coarser mesh $\ThH=\{\cK\}$ which typically
consists of polygonal/polyhedral elements $\cK$ defined as a union of some $K\in\Om_i$,
$i=1,\dots,N$.
We assume that any $\cK\in\ThH$ belongs to only one subdomain $\Om_i$,
including the case $\cK=\oOm_i$ for some $i=1,\dots,N$.
Let $\qK= \min_{K\subset \cK} \pK$, $\cK\in\ThH$, we define the coarse finite element space
(cf. \eqref{Shp})
\begin{align}
  \label{SHP}
  \SHP:= 
   \left\{ v \in \Ldo; \, \restr{v}{\cK} \in P^{\qK}(\cK) \, \forall \cK \in \ThH \right\},
\end{align}
and set
$  H:=\max_{\cK\in\ThH} H_{\cK}$ with $H_{\cK}= \diam(\cK),\ \cK\in\ThH$.

For the purpose of numerical analysis in Section~\ref{sec:anal}, we need also the following
restrictions related to the coarse mesh, cf~\cite[Assumption~18]{Houston-book}.
\begin{assumption}
  \label{ass:coar}
  The coarse elements $\cK\in\ThH$  are star-shaped with respect to an inscribed ball
  that has a radius equivalent to diameter $H_{\cK}$.
  Moreover, we assume that there exists a covering $\ThHs = \{S\}$ of $\ThH$
  consisting of shape-regular simplicial elements  $S\in\ThHs$
  such that for any $\cK\in\ThH$,  there exists  $S_\cK\in\ThHs$ satisfying
  $\cK\subset S_\cK$, $\diam(S_\cK) \leq C_\# H_{\cK}$,
  where $C_{\#}$ is independent of the mesh parameters. Finally, we assume that
  \begin{align}
    \max\nolimits_{\cK\in\ThH} \mathrm{card}\{\cK'\in\ThH:\ \cK'\cap S_{\cK} \not=\emptyset\}
  \end{align}
  is bounded by a positive constant  independent of mesh parameters.
\end{assumption}

\subsection{Local problems}
\label{sec:local}
We introduce the restriction operators $\Ri: \Shp \to \Ship$,  $i=1,\dots,N$ as
$\Ri v_h = v_h|_{\Om_i}$, $i=1,\dots, N$, $v_h\in\Shp$. The corresponding prolongation
operators $\RiT:\Ship\to\Shp$ are given by
$\RiT v_h = v_h$ on $\Om_i$ and
$\RiT v_h = 0$ on $\Om\setminus\Om_i$ for
$v_h\in\Ship$, $i=1,\dots,N$.
Moreover, since $\SHP\subset \Shp$, we define the prolongation operator
$\RNT: \SHP \to \Shp$ as a standard injection from $\SHP$ to $\Shp$.
The restriction operator $\RN: \Shp \to \SHP$  is given as the dual operator
to $\RNT$ with respect to the $L^2(\Om)$-duality.

Furthermore, using \eqref{Ah}, we define the {\em local forms} $\Ahi:\Ship\times\Ship\to \R$ by
\begin{align}
  \label{Ahi}
  \Ahi(u_i, v_i) := \Ah(\RiT u_i, \RiT v_i),\quad u_i,v_i\in\Ship,\quad i=0,\dots, N,
\end{align}
and the {\em local projections} $\tPi:\Shp\to\Ship$ and $\PPi:\Shp\to\Shp$ by
\begin{align}
  \label{Phi}
  \Ahi(\tPi u, v_i) = \Ah(u, \RiT v_i)\quad \forall v_i\in\Ship
\quad \mbox{ and } \quad \PPi:=\RiT\tPi,   \quad i=0,\dots, N,
\end{align}
respectively.
Finally, the {\em two-level additive Schwarz preconditioned operator} reads as follows.
\begin{align}
  \label{Padd}
  \Padd:=\sum\nolimits_{i=0}^N \PPi.
\end{align}

\subsection{Algebraic representation}

Let $\Bh:=\{\vp_{k}\}_{k=1}^\dof$ be the basis
of $\Shp$ ($\dof = \dim\Shp$),
the support of 
each basis function $\vp\in\Bh$ is just one $K\in\Th$. 
We assume that the basis functions are numbered such that
we first number the test functions with support in $\Om_1$,
then the functions with support in $\Om_2$, etc.
%

Then problem \eqref{DG} is equivalent to the linear algebraic system
\begin{align}
  \label{alg1}
  \mA \bku = \bkg,
\end{align}
where $\mA$, $\bku$, and $\bkg$ are algebraic representations
of $\Ah$, $\uh$, and $\lh$, respectively, in the basis $\Bh$.
Similarly, let $\mRi$ and $\mRiT$ denote the algebraic
representations of the operators $\Ri$ and $\RiT$ (cf. Section~\ref{sec:local})
respectively, for $i=0,\dots, N$.
In particular,
the matrices $\mRiT$ for $i=1,\dots, N$ are just extensions of unit matrices
of the corresponding sizes by zero blocks and $\mRNT$ is given by an evaluation of
the basis functions of $\SHHP$ by basis functions of $\Bh$,
we refer, e.g. \cite{AntoniettiALL_JSC14} for details. Obviously, 
$\mRi$ are transposed matrices of $\mRiT$ for $i=0,\dots, N$.

Moreover, the algebraic representation of the local bilinear forms $\Ahi$, the operators
$\tPi$, $\PPi$ 
from \eqref{Ahi}--\eqref{Phi} reads
(cf. \cite{Karashian2001,AntoniettiHouston_JSC11,AntoniettiALL_IJNAM16,DryjaKrzyzanowski_NM16})
\begin{align}
  \label{repre}
  &\mA_i= \mRi\mA\mRiT,\qquad
  \mtPi = \mA_{i}^{-1} \mRi \mA, \qquad
  \mPPi = \mRiT \mA_{i}^{-1} \mRi \mA, \quad  i=0,\dots, N, 
\end{align}
and the representation of the additive Schwarz operator $\Padd$ from \eqref{Padd} is
\begin{align}
  \label{repre2}
  \mPadd=\sum\nolimits_{i=0}^N \mPPi
  = \sum\nolimits_{i=0}^N \mRiT \mA_{i}^{-1} \mRi \mA =: \mNadd^{-1} \mA.
\end{align}
The matrix $\mNadd^{-1} = \mNN^{-1} + \mNad^{-1}$ is the two-level additive Schwarz preconditioner
where 
\begin{align}
  \label{repre3}
  \mNN^{-1}:= \mRNT \mA_{0}^{-1} \mRN
  \quad \mbox{and}\quad
  \mNad^{-1}:=   \sum\nolimits_{i=1}^N \mRiT \mA_{i}^{-1} \mRi = \diag(\mA_1^{-1},\dots, \mA_N^{-1}).
\end{align}
The last matrix-block equality above is valid due to the
numbering of basis functions used. 
Moreover, we set
\begin{align}
  \label{repre4}
  \mPad:=\sum\nolimits_{i=1}^N \mPPi = \sum\nolimits_{i=1}^N \mRiT \mA_{i}^{-1} \mRi \mA
  = \mNad^{-1} \mA.
\end{align}
Finally, we note that $\mPPi$ 
are projections since, by virtue of \eqref{repre},
we have
\begin{align}
  \label{proj2}
  \mPPi^2 =  \mRiT \mA_{i}^{-1} \mRi \mA \,  \mRiT \mA_{i}^{-1} \mRi \mA
  = \mRiT \mA_{i}^{-1}\mA_{i}\mA_{i}^{-1} \mRi \mA
  = \mPPi, \quad i=0,\dots,N
\end{align}
provided that the performance of local solvers $\mA_{i}^{-1}$ is carried out exactly.



\subsection{Iterative methods}

To better understand the performance of the preconditioner, we consider
the (one- and two-level) iterative Schwarz schemes to solve \eqref{alg1}.
The matrix $\mA$ has blocks $\mA=\{\mA_{ij}\}_{i,j=1}^N$
such that each block $\mA_{ij}$ contains entries corresponding to $\Ah(\vp_{k},\vp_{l})$,
$k,l=1,\dots,\dof$ satisfying $\supp(\vp_{k}) \subset \oOm_j$
and $\supp(\vp_{l}) \subset \oOm_i$.  
Obviously, if $i\not=j$ and 
$\Om_i$ and $\Om_j$ have no common edge ($d=2$) or face ($d=3$)
then $\mA_{ij}=\bfzero$. Moreover, due to \eqref{repre},
the diagonal blocks $\mA_{ii} = \mA_i$, $i=1,\dots,N$, and $\mA_{ij} = \mRi\mA\mRjT$, $i,j=1,\dots,N$.


First, we mention the {\em one-level additive Schwarz method} that solves \eqref{alg1} iteratively.
Let $\bku^\ell = (\bku_1^{\ell},\dots, \bku_N^\ell)^\T$ denote the $\ell$ -th approximation of
$\bku$, then we have 
\begin{align}
  \label{dd2}
  \bku^{\ell+1} = \bku^{\ell} + \mNad^{-1}( \bkg   -  \mA\bku^{\ell}),
  \qquad \ell=0,1,2,\dots,
\end{align}
where $\mNad^{-1}$ is given by \eqref{repre3}.
Obviously, due to discontinuous Galerkin discretization,
scheme \eqref{dd2} is equivalent to the block Jacobi method. 

Moreover, the {\em two-level additive Schwarz method} reads as follows: For $\ell=0,1,2,\dots$, set
\begin{align}
  \label{2add4}
  \left.
  \begin{array}{rl}
   \mbox{(i)}   & \quad \bku^{\ell+1/2}:= \bku^{\ell} + \mNN^{-1} (\bkg - \mA \bku^{\ell}) \\
   \mbox{(ii)}  & \quad\bku^{\ell+1} := \bku^{\ell+1/2} + \mNad^{-1}( \bkg   -  \mA\bku^{\ell}) \\
  \end{array}
  \right\}
  \ \Longleftrightarrow \ 
  \bku^{\ell+1} = \bku^{\ell} + \mNadd^{-1}( \bkg   -  \mA\bku^{\ell})
\end{align}
due to \eqref{repre2} -- \eqref{repre3}.
Obviously, if $\bku^\ell \to \bku$ for $\ell \to \infty$, then the limit $\bku$
is also the solution of \eqref{alg1}. Moreover, with respect to \eqref{repre2} and \eqref{2add4},
the vector $\bku$
satisfies
\begin{align}
  \label{2add5}
  \mPadd \bku = \mNadd^{-1} \mA \bku = \mNadd^{-1} \bkg.
\end{align}
Therefore, the two-level additive Schwarz preconditioner \eqref{repre2} corresponds to
the iterative method \eqref{2add4}.  
We note that method \eqref{2add4} does not converge
if the spectral radius satisfies $\rho(\mNadd^{-1} \mA) > 2$. 
This may occur in the general case; therefore, the matrix $\mNadd^{-1}$ is typically used
as a preconditioner for an iterative Krylov solver.
For further details, we refer the reader to \cite{bjorstad1991spectra}.

Finally, we introduce the {\em hybrid Schwarz method} which was introduced 
in \cite{Mandel_94} in the context of the conforming finite element method, see also
\cite[Section~2.5.2]{ToselliWidlund-DD05}. It can be written as a three-sub-step scheme
as follows.
\begin{subequations}
  \label{hybM1}
  \begin{align}
    \label{hybM1a}
    \mbox{(i)} & \quad
    \bku^{\ell+1/3} := \bku^{\ell} + \mNN^{-1} (\bkg - \mA \bku^{\ell}), \\
    \label{hybM1b}
    \mbox{(ii)} & \quad
    \bku^{\ell+2/3} := \bku^{\ell+1/3} + \mNad^{-1}( \bkg   -  \mA\bku^{\ell+1/3}), 
    \\
    \label{hybM1c}
    \mbox{(iii)} & \quad
    \bku^{\ell+1} := \bku^{\ell+2/3} + \mNN^{-1} (\bkg - \mA \bku^{\ell+2/3}), \hspace{12mm}
    \ell =0,1,2, \dots.
  \end{align}
\end{subequations}
The first and third steps correspond to the application of the global coarse
operator, and the middle step to the local ones. This performance is symmetric, and therefore
it can be applied as the preconditioner for symmetric problems. We note that whereas
both steps in the additive method \eqref{2add4} can be carried out simultaneously,
the hybrid method \eqref{hybM1}
does not allow us to solve the coarse problem in parallel with the local ones.
However, in order to achieve a practical scalability of the computation, the global problem
is chosen to be typically smaller than the local ones.

Now, we derive the explicit form of the hybrid preconditioner.
Subtracting the exact solution $\bku$ from \eqref{hybM1a} and using \eqref{alg1}, we have
\begin{align}
  \label{hybM2}
  \bku^{\ell+1/3} -\bku = \bku^{\ell} -\bku + \mNN^{-1} (\mA \bku - \mA \bku^{\ell}) 
  \quad \Longrightarrow \quad
  \bke^{\ell+1/3}  = (\mI - \mNN^{-1} \mA) \bke^{\ell},
\end{align}
where $\mI$ is the unit matrix and
$\bke^j := \bku^j - \bku$, $j\in\{\ell,\ell+1/3,\ell+2/3,\ \ell=0,1,\dots\}$ denotes
the error of the $j$-th approximation.
A similar procedure for steps \eqref{hybM1b} and \eqref{hybM1c}, together
with \eqref{repre3} -- \eqref{repre4}, gives the expression for the propagation of errors in one step.
\begin{align}
  \label{hybM3}
  \bke^{\ell+1}  = (\mI - \mNN^{-1} \mA)  (\mI - \mNad^{-1} \mA)  (\mI - \mNN^{-1} \mA) \bke^{\ell}
   = (\mI - \mPPN)  (\mI - \mPad)  (\mI - \mPPN) \bke^{\ell}.
\end{align}
Formula \eqref{hybM3} can be written as 
\begin{align}
  \label{hybM4}
  \bke^{\ell+1}  = (\mI - \mPhyb) \bke^{\ell},\qquad \mbox{where}\quad
  \mPhyb := \mI - (\mI - \mPPN)  (\mI - \mPad)  (\mI - \mPPN),
\end{align}
and implies (similarly as in \eqref{hybM2})
\begin{align}
  \label{hybM5}
  \bku^{\ell+1} - \bku  = \bku^{\ell} - \bku  + \mPhyb (\bku - \bku^{\ell})
  \quad \Rightarrow \quad
  \bku^{\ell+1} = \bku^{\ell}  + \mNhyb^{-1} (\bkg - \mA \bku^{\ell}).
\end{align}
The matrix $\mNhyb^{-1}$ is the hybrid Schwarz preconditioner given by
$\mPhyb = \mNhyb^{-1} \mA$.
Using \eqref{repre} -- \eqref{repre4}, we can evaluate
the explicit formulas for $\mPhyb$ and $\mNhyb^{-1}$ as
 \begin{align}
   \label{hyb3}
   \mPhyb & = \mPPN +  \big(\mI - \mPPN \big) \mPad \big(\mI - \mPPN \big)
   =  \mNN^{-1}\mA +  \big(\mI -  \mNN^{-1}\mA \big)  \mNad^{-1}\mA \big(\mI - \mNN^{-1}\mA \big),
 \end{align}
 and
 \begin{align}
   \label{hyb4}
    \mNhyb^{-1} & =
   \mNN^{-1} + \big(\mI - \mNN^{-1} \mA \big) \mNad  \big(\mI - \mA \mNN^{-1} \big)   \\
   & = \mRNT \mA_{0}^{-1} \mRN + \big(\mI - \mRNT \mA_{0}^{-1} \mRN\mA \big)
  \sum\nolimits_{i=1}^N \mRiT \mA_{i}^{-1} \mRi \big(\mI - \mA \mRNT \mA_{0}^{-1} \mRN\big). \notag
\end{align}

Algorithms~\ref{alg:ASM} and \ref{alg:ASM3} describe
the applications of preconditioners $\mNadd^{-1}$ and $\mNhyb^{-1}$ introduced in
\eqref{repre2} and \eqref{hyb4}, respectively.
As mentioned above,
Algorithm~\ref{alg:ASM} allows one to solve local fine problems
(with $\mA_i$, $i=1,\dots,N$) together with the global coarse problem (with $\mA_0$)
in parallel.
On the other hand,
Algorithm~\ref{alg:ASM3} requires first solving the global coarse problem,
then the local fine problems, and finally
the global coarse problem. In addition, two multiplications by
$\mA$ have to be performed. Multiplication by $\mA$ is typically much cheaper than
solving local or global problems. Moreover, if the global coarse problem is smaller
than the local ones, then the application of the preconditioner 
$\mNhyb^{-1}$
exhibits only a small increase in computational time compared to
the preconditioner $\mNadd^{-1}$. The study of the computational costs of the
additive and hybrid preconditioner is given in Section~\ref{sec:numer}.
  
\begin{algorithm}
  \caption{Application of preconditioner $\mNadd^{-1}$ from \eqref{repre2}:
    $\bku \leftarrow \mNadd^{-1}\bkx$}
  \label{alg:ASM}
  \begin{spacing}{1.15}
    \begin{algorithmic}[1]
      \STATE \textbf{input} matrices $\mA$, $\mA_{i}$, $\mRi$, $i=0,\dots,N$, vector $\bkx$
      \STATE $\bkx_i := \mRi \bkx$ and solve $\mA_{i} \bky_i = \bkx_i$ for $i=0,\dots,N$
      \STATE \textbf{output} vector  $\bku := \sum_{i=0}^N \mRi^T \bky_i$
    \end{algorithmic}
  \end{spacing}
\end{algorithm}


\begin{algorithm}
  \caption{Application of preconditioner $\mNhyb^{-1}$ from \eqref{hyb4}:
    $\bku \leftarrow \mNhyb^{-1}\bkx$}
  \label{alg:ASM3}
  \begin{spacing}{1.15}
    \begin{algorithmic}[1]
      \STATE \textbf{input} matrices $\mA$, $\mA_{i}$, $\mRi$, $i=0,\dots,N$, vector $\bkx$
      \STATE $\bkx_0 := \mRN \bkx$  and solve $\mA_{0} \bky_0 = \bkx_0$ \label{asm3_c1}
      \STATE $\bkz_0 :=  \mRN^T \bky_0$,  
            $\bkz := \bkx - \mA \bkz_0$  \label{asm3_m2}
      \STATE $\bkz_i := \mRi \bkz$  and solve $\mA_{i} \bky_i = \bkz_i$ for $i=1,\dots,N$
      \STATE $\bky := \sum_{i=1}^N \mRi^T \bky_i$
      \STATE $\bkw_0 : = \mRN \mA \bky$ and solve $\mA_{0} \bkv_0 = \bkw_0$ \label{asm3_c2}
      \STATE \textbf{output} vector $\bku = \bkz_0 + \bky - \mRN^T \bkv_0$
    \end{algorithmic}
  \end{spacing}
\end{algorithm}

\section{Numerical analysis}
\label{sec:anal}

In this section, we derive the spectral bound of the hybrid Schwarz operator \eqref{hyb3}.
The notation $a \lesssim b$ means that there exists a constant $C$ independent of the
discretization parameters such that $a \leq C\, b$.
Moreover, $a \approx b$ means $a \lesssim b$ and $b \lesssim a$.
For the sake of simplicity, we consider
quasiuniform meshes $\Th$ and $\ThH$ with mesh steps $h$ and $H$, respectively,
and constant polynomial degrees $p\ge q\ge 1$. 
Hence, we have $\pg\approx p$ and $\hg\approx h$ for all $\gamma \in\Fh$.
Moreover, we recall the symbols for the upper and lower bounds of the diffusion
$\barK$ and $\ubarK$, respectively, given by \eqref{Kij}.


\subsection{Auxiliary results}

By symbol $\nabla_h v_h$, we denote the piecewise element gradient operator for $v_h\in\Shp$.
For the function $v_h\in\Shp$, we define the DG-norm defined as
  \begin{align}
  \label{DGnorm}
    \dg \vh \dg^2 = 
  \ubarK\bnormP{ \nabla_h \vh}{\Ldo}{2}
  + \sumfh   \bnormP{\sqrt{\sigma} \ \jump{\vh}}{\Ldd}{2}
  \ \mbox{ with } \ 
  \restr{\sigma}{\gamma} = \CW \K_\gamma {p^2}/{h},
  \end{align}
  cf.~\eqref{sigma}.
  The form $\Ah$ from \eqref{Ah} is coercive and continuous \cite[Chapter~2]{DGM-book}, namely
\begin{align}
  \label{coerc}
  \dg \uh \dg^2 & \lesssim \Ah(\uh,\uh),\qquad \uh\in\Shp, \\
  \label{cont1}
  |\Ah (\uh,\vh)| & \lesssim \frac{\barK}{\ubarK}
  \dg \uh \dg\, \dg \vh \dg,\qquad \uh,\vh \in \Shp, 
\end{align}

\begin{remark}
  \label{rem:K_const}
  We note that if the diffusion function $\K$ is piecewise constant, 
  we can replace the DG-norm \eqref{DGnorm} by
  \begin{align}
  \label{DGnorm0}
  \dg \vh \dg^2 = 
  \bnormP{\sqrt{\K} \ \nabla_h \vh}{\Ldo}{2}
  + \sumfh   \bnormP{\sqrt{\sigma} \ \jump{\vh}}{\Ldd}{2}
  \ \mbox{ with } \ 
  \restr{\sigma}{\gamma} = \CW {\K}_\gamma {p^2}/{h}.
\end{align}
Then the  coercivity and continuity results read
\begin{align}
  \label{coerc0}
  \dg \uh \dg^2 & \lesssim \Ah(\uh,\uh),\qquad \uh\in\Shp, \\
  \label{cont0}
  |\Ah (\uh,\vh)| & \lesssim  
  \dg \uh \dg\, \dg \vh \dg,\qquad \uh,\vh \in \Shp, 
\end{align}

\end{remark}

Moreover, we consider a ``diffusion-independent'' norm
\begin{align}
  \label{DGnorm2}
  \dg \vh \dgI^2 :=
  \bnormP{ \nabla_h \vh}{\Ldo}{2}
  + \sumfh   \bnormP{\sqrt{\sbar} \ \jump{\vh}}{\Ldd}{2}
  \ \mbox{ with } \
  \restr{\sbar}{\gamma} = \CW {p^2}/{h}.
\end{align}
Obviously, due \eqref{Kij},
\begin{align}
  \label{DGnorm2a}
  \dg \vh \dg^2 \leq \barK \dg\vh\dgI^2,\qquad
  \dg\vh\dgI^2 \leq \frac{1}{\ubarK} \dg \vh \dg^2,
  \qquad
  \vh\in\Shp
\end{align}
where $\dg\cdot\dg$ is given by \eqref{DGnorm} or by \eqref{DGnorm0} if $\K$ is piecewise
constant.

We recall the following variants of the trace inequalities for piecewise polynomial
functions.
\begin{lemma}
  \label{trace_ineq}
  Let  meshes $\Th$ be shape-regular and the coarse mesh $\ThH$
  fulfill Assumption~\ref{ass:coar}.
  For any $\vh \in \Shp$, we have 
  \begin{align}
    \label{inverse}
    \normP{\vh}{\LdK}{2} & \lesssim  \frac{p^2}{h} \normP{\vh}{L^2(K)}{2},
    \qquad K\in\Th,\ \vh\in\Shp, \\
    \label{trace1}
    \sum_{\cK \in \ThH} \normP{\vh}{\LdcK}{2}& \lesssim    \norm{\nabla_h\vh}{\Ldo} \norm{\vh}{\Ldo}
    + \frac{1}{H} \normP{\vh}{\Ldo}{2} \\
    & +
    \bigg(\sum_{\cK \in \ThH} \sum_{\gamma\in\Fh,\gamma\subset\cK}  \bnormP{\sbar^{1/2}\jump{\vh}}{\Ldd}{2}
    \bigg)^{1/2}  \norm{\vh}{\Ldo} , \quad \vh\in\Shp. \notag
  \end{align}
\end{lemma}
{\it Proof:}\ 
  Inequality \eqref{inverse} is the standard inverse / trace inequality 
  (e.g. \cite{riviere-book,Schwab-book}).
  The estimate \eqref{trace1} is proved in \cite[Lemma~5]{Smears_JSC18}.\qed

Moreover, we employ the following result concerning the approximation of the function of $\Shp$
by a function of $\SHP$, cf.~\eqref{SHP}.

\begin{lemma}{\cite[Lemma~5.1]{AntoniettiALL_IJNAM16}}
    \label{lem:bound}
    Let  meshes $\Th$ be shape-regular and the coarse mesh $\ThH$
    fulfill Assumption~\ref{ass:coar}.
    Let $\wh \in \Shp$, then there exists $\wH \in \SHP$ such that the following inequalities holds
    \begin{align}
      \label{boundL2}
      \norm{\wh - \wH}{\Ldo} & \lesssim \frac{H}{q} \dg \wh \dgI \\
      \label{boundH1}
      \norm{\nabla_h(\wh - \wH)}{\Ldo} & \lesssim \dg \wh \dgI. 
    \end{align} 
\end{lemma}

{\it Proof:}\ 
  The proof is based on an approximation of $\wh \in \Shp$ by function $w\in H^1(\Om)$, and
  then $\wH$ is defined as the $L^2(\Om)$ projection of $w$ in $\SHP$. For details,
  we refer to \cite[Lemma~5.1]{AntoniettiALL_IJNAM16}.
  See also \cite[Lemma~5.13]{AntoniettiALL_MC20},
  where the local variants of these estimates related to the subdomains $\Om_i$,
  $i=1,\dots,N$ are presented.\qed

Consequently, we have the following auxiliary estimate.
\begin{lemma}
  \label{lem:L2d}
  Let $\wh \in \Shp$ and $\wH \in \SHP$ be given by Lemma~\ref{lem:bound}. Then
  \begin{align}
    \label{eq:L2d}
    \sum_{\cK \in \ThH}\bnormP{\sigma^{1/2}(\wh-\wH)}{\LdcK}{2}
    \lesssim  \frac{\barK}{\ubarK} \frac{p^2}{q} \frac{H}{h} \Ah(\wh,\wh). 
  \end{align}
\end{lemma}
{\it Proof:}\ 
  follows from the trace inequality \eqref{trace1}.
  Let $\wh \in \Shp$ and $\wH \in \SHP$ be given by Lemma~\ref{lem:bound}.
  Then $\jump{\wH}_\gamma =0$ for $\gamma\in\Fh$, $\gamma\subset\cK$, $\cK\in\ThH$.
  Using this fact and \eqref{DGnorm2}, we bound the ``face'' term in \eqref{trace1} for
  $\vh:=\wh-\wH$ as
  \begin{align}
    \label{jump1}
    \sum_{\cK \in \ThH} \sum_{\gamma\in\Fh,\gamma\subset\cK}
    \bnormP{\sbar^{1/2}\jump{\wh-\wH}}{\Ldd}{2}
    \leq
    \sum_{\gamma\in\Fh}
    \bnormP{\sbar^{1/2}\jump{\wh}}{\Ldd}{2} \leq  \dg \wh \dgI^2.
  \end{align}
  Moreover,  putting $\vh:=\wh-\wH$ in  \eqref{trace1}, using \eqref{DGnorm2},
  \eqref{boundL2} -- \eqref{boundH1},
  \eqref{DGnorm2a} and \eqref{jump1}, we obtain
  \begin{align}
    \label{jump2}
    &\sum_{\cK \in \ThH}\bnormP{\sigma^{1/2}(\wh-\wH)}{\LdcK}{2} \\
    &
    \lesssim
        {\barK}\frac{p^2}{h} \norm{\wh-\wH}{\Ldo} \left(
        \norm{\nabla_h(\wh-\wH)}{\Ldo} 
        + \frac{1}{H} \normP{\wh-\wH}{\Ldo}{}
        +  \dg \wh \dgI
        \right) \notag\\
        &
        \lesssim
            {\barK} \frac{p^2}{h}
            \frac{H}{q} \left( 1 + \frac{H}{Hq}  + 1\right)\dg \wh \dgI^2
            \lesssim \frac{\barK}{\ubarK}   \frac{p^2}{h} \frac{H}{q}  \dg \wh \dg^2,
             \notag
  \end{align}
 which together with the coercivity \eqref{coerc} proves lemma. \qed

The following result is a small modification of \cite[Lemma~4.2]{AntoniettiHouston_JSC11},
see also \cite[Lemma~5.14]{AntoniettiALL_MC20}.

\begin{lemma}
    \label{lem:decom}
    Let $v \in \Shp$, it has a unique decomposition $v = \sum_{i=1}^{N} \RiT \vi$
    with $\vi \in \Ship$. Then 
    \begin{align}
      \label{st}
      \bigg\lvert \sum_{\substack{i,j = 1\\ i \neq j}}^N \Ah(\RiT \vi, \RjT \vj) \bigg\rvert
      \lesssim 
      \barK   \normP{\nabla_h v}{L^2(\Om)}{2}
      + \sum_{\cK \in \ThH}\bnormP{\sigma^{1/2}v}{\LdcK}{2},
    \end{align}
    where $\barK$ is given by \eqref{Kij}.
\end{lemma}
{\it Proof:}\ 
  Let $v = \sum_{i=1}^{N} \tvi$, where $\tvi :=\RiT \vi$, $i=1,\dots,N$ is
  the decomposition of $v\in\Shp$.
  From \eqref{Ah}, we observe that $\Ah(\tvi, \tvj)$ is non-vanishing for $i\not=j$,
  only if $\Om_i$ and $\Om_j$ have a common boundary $\gomij:=\gom_i\cap\gom_j$
  (with positive $(d-1)$-dimensional measure).
  In this case, we have
  \begin{align}
    \label{st2}
    \Ah(\tvi,\tvj) = &
    - \sumFij \Lspg{ \aver{\K\nabla \tvi}}{\jump{\tvj}}
    - \sumFij\Lspg{ \aver{\K\nabla \tvj} }{\jump{\tvi}}
    +\sumFij \Lspg{\sigma\jump{\tvi}}{\jump{\tvj}}.
  \end{align}
  We bound the first term in \eqref{st2}. The Cauchy-Schwarz and Young inequalities imply
  \begin{align}
    \label{st3}
    \Big |\sumFij \Lspg{ \aver{\K\nabla \tvi}}{\jump{\tvj}} \Big|
    \lesssim
    \bnormP{\sigma^{-1/2}  \K\nabla \tvi }{L^2(\gom_i)}{2}
    + \bnormP{\sigma^{1/2} \tvj}{L^2(\gom_j)}{2},
  \end{align}
  since $\tvi=0$ on $\gom_j$ and $\tvj=0$ on $\gom_i$. Moreover, due to  $\K_\gamma \ge \K$
  on $\gom_i$, and \eqref{inverse}, we have
  \begin{align}
    \label{st4}
    \normP{\sigma^{-1/2}  \K\nabla \tvi }{L^2(\gom_i)}{2} 
    & \approx \sum_{\gamma\subset \gom_i} \int_{\gamma} \frac{h}{p^2 \K_\gamma } \K^2 |\nabla \tvi|^2\dS 
    \lesssim \sum_{\gamma\subset \gom_i} \int_{\gamma} \frac{h}{p^2 } \K |\nabla \tvi|^2\dS \\
    & \leq \barK \sum_{K\in\Thi} \frac{h}{p^2} \normP{\nabla \tvi}{L^2(\dK)}{2}
    \lesssim \barK \sum_{K\in\Thi} \normP{\nabla \tvi}{L^2(K)}{2}, \notag
  \end{align}
  which together with \eqref{st3} gives
    \begin{align}
    \label{st5}
    \Big |\sumFij \Lspg{ \aver{\K\nabla \tvi}}{\jump{\tvj}} \Big|
    \lesssim \barK  \normP{\nabla \tvi}{L^2(\Om_i)}{2}
    + \bnormP{\sigma^{1/2} \tvj}{L^2(\gom_j)}{2}.
    \end{align}
    In the same manner as \eqref{st5}, we estimate the second term in \eqref{st2}.
    Additionally, applying the Cauchy and Young inequalities to the third term of \eqref{st2}, we have
    \begin{align}
      \label{st6}
      \Big|\sumFij \Lspg{\sigma\jump{\tvi}}{\jump{\tvj}} \Big|
      \lesssim \bnormP{\sigma^{1/2} \tvi}{L^2(\gom_i)}{2}
      + \bnormP{\sigma^{1/2} \tvj}{L^2(\gom_j)}{2},
    \end{align}
    which together with \eqref{st5} gives \eqref{st}.
  \qed

\subsection{Main theoretical results} We derive a bound of the preconditioned operator
using the modification of the approach from \cite[Chapter~2]{ToselliWidlund-DD05}.
First, the following strengthened Cauchy-Schwarz inequalities which
simply follow from the properties of the local spaces
$\Ship$, $i=1,\dots,N$.

\begin{lemma}
  \label{assum2}
    There exist constants $0 \leq \ve_{ij} \leq 1$, $i,j = 1, \dots ,N$, such that 
    \begin{align}
      \label{CS}
      \vert \Ah (\RiT u_i,\RjT u_j) \vert \leq \ve_{ij}
      \Ah (\RiT u_i,\RiT u_i)^{1/2} \Ah (\RjT u_j,\RjT u_j)^{1/2}, \quad i,j = 1, \dots N,
    \end{align}
    for all $u_i \in \Ship$, $u_j \in \Shjp$. Moreover, let $\rho(\bve)$ be the spectral radius of    $\bve = \{ \ve_{ij} \}_{i,j = 1}^N$ then
    \begin{align}
      \label{CS1}
      \rho(\bve)= \NS + 1,
    \end{align}
    where $\NS$ is the maximum number of adjacent subdomains of any given subdomain in domain decomposition.
\end{lemma}
{\it Proof:}\ 
  Since $\Ah$ defines a scalar product, then $\ve_{ij}=1$ for either $i=j$ or if $\Om_i$ and $\Om_j$ share a boundary.
  Otherwise, $\ve_{ij}=0$. \qed

Furthermore, we present the stability of the decomposition.
\begin{lemma}
  \label{lem3}
  There exists a constant $\Co$, such that
  for  any  $u\in\Shp$,  
  there exists the decomposition
  $u = \sum_{i=0}^{N} \RiT u_i, \, u_i \in \Ship$ that satisfies
  \begin{align}
    \label{eq:assum3}
    \sum\nolimits_{i=0}^{N} \Ahi(u_i,u_i) \leq \Co \Ah (u,u),
  \end{align}
  with the constant
  \begin{align}
    \label{C0}
    \Co = C_\sigma \frac{\barK}{\ubarK}  \frac{H}{h} \frac{p^2}{q},
  \end{align}
  where
  $C_\sigma>0$ is independent of the mesh sizes $h$, $H$ and the polynomial degrees $p$, $q$.
\end{lemma}

{\it Proof:}\ 
  The proof follows closely the proof presented in \cite[Theorem~5.1]{AntoniettiALL_IJNAM16}
  where the modified versions of the auxiliary lemmas (cf. Lemmas~\ref{lem:L2d} and \ref{lem:decom})
  are employed. Therefore, we skip the proof and refer to \cite{AntoniettiALL_IJNAM16}.
  \qed

We derive the bound of the hybrid preconditioner.

\begin{theorem}
  \label{thm:hyb}
  Let $\Phyb$ be the hybrid preconditioned operator (cf.~\eqref{hyb3})
  \begin{align}
  \label{HYB}
  \Phyb = \PPN +  \big(I - \PPN \big)\sum\nolimits_{i=1}^N \PPi \big(I - \PPN \big).
\end{align}
  Then
  \begin{align}
    \label{SB}
    \max \{1, \Coo^{2} \}^{-1} \Ah(u,u) \leq \Ah(\Phyb u, u) \leq \max \{1, \rho(\bve) \} \Ah(u,u)\qquad \forall u\in\Shp,
  \end{align}
  where $\rho(\bve)$ and  $\Coo$ are discretization parameters given by \eqref{CS1} and
  \eqref{C0}, respectively.
\end{theorem}

\begin{remark}
  \label{rem:Widlund}
  The assertion \eqref{SB} is the same as \cite[Theorem~2.13]{ToselliWidlund-DD05}.
  However, it is based on a different assumption (\cite[Theorem~2.12]{ToselliWidlund-DD05}),
  namely, the existence of
  $\Coo>0$ such that any $u\in  \range(I-\PPN)$ admits a decomposition
  $u = \sum_{i=1}^N \u_i$ that fulfills $\sum\nolimits_{i=1}^{N} \Ahi(u_i,u_i) \leq \Co \Ah (u,u)$,
  compare with cf.~\eqref{eq:assum3}. We are not able to verify the original assumption for the
discontinuous Galerkin discretization used. Therefore, we use the decomposition with estimate~\eqref{eq:assum3}
and modify the proof of \cite[Theorem~2.13]{ToselliWidlund-DD05}.
  Although it is very similar, we present it for completeness.
\end{remark}

{\it Proof:}\ [of Theorem~\ref{thm:hyb}]
  Let $v\in\Shp$, from \eqref{Ahi} -- \eqref{Phi}, we obtain the identity
  \begin{align}
    \label{SB1}
    \Ah(\PPi v,\PPi v ) = \Ah(\RiT\tPi v,\RiT\tPi v ) = \Ahi(\tPi v,\tPi v )
    = \Ah(v, \RiT\tPi v ) = \Ah(v, \PPi v )
  \end{align}
  for $i=1,\dots,N$.
  Let $\PPh = \sum_{i=1}^N\PPi$.  From the strengthened Cauchy-Schwarz inequalities
  \eqref{CS} -- \eqref{CS1} 
  and \eqref{SB1}, we have
  \begin{align}
    \label{SB2}
    \Ah(\PPh v,\PPh v) & = \sum_{i,j=1}^N \Ah(\PPi v, \PPj v)
    \leq \sum_{i,j=1}^N \ve_{ij} \Ah(\PPi v, \PPi v)^{1/2} \Ah(\PPj v, \PPj v)^{1/2} \\
    & \leq \rho(\bve) \sum_{i=1}^N  \Ah(\PPi v, \PPi v)
    = \rho(\bve) \sum_{i=1}^N  \Ah( v, \PPi v) 
    = \rho(\bve)  \Ah( v, \PPh v)  \notag \\
    &\leq  \rho(\bve)  \Ah( v, v)^{1/2} \Ah(\PPh v, \PPh v)^{1/2}. \notag
  \end{align}
  Therefore, using \eqref{SB2}, we deduce
  \begin{align}
    \label{SB3}
    \Ah(\PPh v, v) \leq  \Ah(\PPh v, \PPh v)^{1/2}  \Ah( v, v)^{1/2}\leq \rho(\bve)  \Ah( v, v).
  \end{align}
  Let $u\in\Shp$, we insert $v:=(I-\PPN)u$ in \eqref{SB3}, and since $\PPN$ is self-adjoint,
  we have
  \begin{align}
    \label{SB4}
    \Ah((I-\PPN) \PPh (I-\PPN)u, u) = 
    \Ah(\PPh (I-\PPN)u, (I-\PPN)u) 
    \leq \rho(\bve)  \Ah( (I-\PPN)u, (I-\PPN)u). 
  \end{align}
  Furthermore, from \eqref{hyb3}, the fact that $\PPN$ is orthogonal projection,
  and \eqref{SB4}, we obtain
  \begin{align}
    \label{SB5}
    \Ah(\Phyb u, u) & = \Ah(\PPN u, u) + \Ah((I-\PPN) \PPh (I-\PPN)u, u)  \\
    & \leq \Ah(\PPN u, \PPN u) + \rho(\bve)  \Ah( (I-\PPN)u, (I-\PPN)u) \leq \max \{ 1, \rho(\bve) \} \Ah(u,u),
    \notag
  \end{align}
  which proves the second inequality in \eqref{SB}.
  \par
  Moreover, let $v\in\Shp$ and  $v = \sum_{i=0}^{N} \RiT v_i, \, v_i \in \Ship$ be the decomposition
  from  Lemma~\ref{lem3}. Then from \eqref{Ahi}, \eqref{Phi}, and
  \eqref{eq:assum3}, we derive
  \begin{align}
    \label{SB10}
    \Ah(v,v) &= \sum_{i=0}^{N}\Ah(v, \RiT v_i) =  \sum_{i=0}^{N}\Ahi(\tPi v, v_i)
    \leq \bigg(\sum_{i=0}^{N} \Ahi(\tPi v, \tPi v)\bigg)^{1/2}
    \bigg(\sum_{i=0}^{N} \Ahi(v_i, v_i)\bigg)^{1/2} \notag  \\ 
    & \leq \bigg(\sum_{i=0}^{N} \Ahi(\tPi v, \tPi v)\bigg)^{1/2} \Coo \Ah(v, v)^{1/2}
    = \Coo \Ah(v, v)^{1/2} \bigg( \sum_{i=0}^{N} \Ah( v, \RiT\tPi v)\bigg)^{1/2}  \notag \\ 
    & = \Coo \Ah(v, v)^{1/2} \bigg(\sum_{i=0}^{N} \Ah( v, \PPi v)\bigg)^{1/2}.
  \end{align}
  Obviously, for $v\in \range(I-\PPN)$,  $\PPN v = 0$ and the term for $i=0$ in the last sum is
  vanishing. Therefore, the squaring  estimate \eqref{SB10} gives
 \begin{align}
    \label{SB10a}
    \Ah(v,v) \leq \Coo^2\sum_{i=1}^{N} \Ah( v, \PPi v)
    =  \Coo^2 \Ah( v, \PPh v)
    \qquad
    \mbox{for }v\in \range(I-\PPN).
 \end{align}
 Again, let $u\in\Shp$, we insert $v:=(I-\PPN)u$ in \eqref{SB10a} and obtain
  \begin{align}
    \label{SB11}
    \Ah(\PPh (I-\PPN)u, (I-\PPN)u) \geq \Coo^{-2} \Ah((I-\PPN)u, (I-\PPN)u).
  \end{align}
  Then, similarly as in \eqref{SB5},
   \begin{align}
     \label{SB12}
     \max \{1, \Co \}  \Ah(\Phyb u, u) & \ge  \Ah(\PPN u,\PPN u) + \Co \Ah((I-\PPN) \PPh (I-\PPN)u, u) \\
     & \ge  \Ah(\PPN u,\PPN u) +  \Ah((I-\PPN) u,  (I-\PPN)u) = \Ah(u,u), \notag
   \end{align}
   which proves the first inequality in \eqref{SB}. \qed

\begin{corollary}
  Theorem~\ref{thm:hyb} with  Lemmas~\ref{assum2} and \ref{lem3} implies that 
  the condition number of the hybrid Schwarz operator \eqref{HYB} satisfies
  \begin{align}
    \label{RES1}
    \kappa(\Phyb) :=
    {\max_{v\in\Shp, v\not = 0} \frac{\Ah(\Phyb v, v)}{\Ah(v,v)}} \, /\, 
    {\min_{v\in\Shp, v\not = 0} \frac{\Ah(\Phyb v, v)}{\Ah(v,v)}}
      \lesssim \frac{\barK}{\ubarK} \frac{H}{h}\,\frac{p^2}{q}\, (\NS + 1),
  \end{align}
  where $h$, $p$, $H$, $q$ are the discretization parameters and
  $\NS$ is the maximum number of neighboring subdomains.
\end{corollary}

\begin{remark}
  The spectral bound of the additive Schwarz operator \eqref{Padd} reads
  \begin{align}
    \label{RES2}
    \kappa(\Padd) \lesssim \frac{\barK}{\ubarK} \frac{H}{h}\,\frac{p^2}{q}\, (\NS + 2),
  \end{align}
  see, e.g.,~\cite{AntoniettiALL_IJNAM16} for the case of $\barK=\ubarK$.
  In practice, the constant $\NS$ is usually between 7 and 12.
  Hence, the bound \eqref{RES1} of the hybrid preconditioner is only slightly
  better than the bound \eqref{RES2} of the additive one.
  They are the same in terms of $h$, $H$, $p$, and $q$.
\end{remark}

\begin{remark}
  \label{rem:materialDD}
  For the massively parallel non-overlapping additive Schwarz method
  (cf.~\cite{AntoniettiALL_MC20}),
  the sharper estimate in terms of diffusion $\K(x)$ was derived.
  Instead of the ratio ${\barK}/{\ubarK}$, the bound contains the term
  \begin{align}
    \label{RES3}
    \max_{\cK\in\ThH} \frac{\max_{x\in \cK} \K(x)}{\min_{x\in \cK} \K(x)}.
  \end{align}
  It is valid under assumption that the coarse elements $\cK\in\ThH$ are convex,
  which can be achieved typically for non-nested $\Th$ and $\ThH$ (which is not our case).
  This approach can give a better spectral bound than \eqref{RES2}, e.g.,
  if the domain decomposition is aligned with the material interfaces.
  Although our theoretical result does not cover this aspect, we present
  its numerical study in Section~\ref{sec:numer}.
\end{remark}

\begin{remark}
 The relationship between the condition numbers of the additive and hybrid Schwarz preconditioners
has also been studied in \cite{Mandel_94}, where a similar approach based on stable decomposition is used.
In that work, the authors derive the following estimate:
  \begin{align}
    \kappa(\Phyb) \leq \kappa(\Padd).
  \end{align}
\end{remark}

\subsection{Computational costs}

Finally, we discuss the computational costs of an iterative solver employing the presented
preconditioners. 
The most expensive part
of Algorithms~\ref{alg:ASM}--\ref{alg:ASM3} is the performance of the preconditioner,
which exhibits the solution of the local (fine) systems and the global (coarse) one
\begin{align}
  \label{flops0}
  \mA_{i} \bky_i = \bkx_i,\quad i=1,\dots,N,\qquad \mbox{and} \qquad
  \mA_{0} \bky_0 = \bkx_0,
\end{align}
respectively. These systems are usually solved by a direct method,
e.g., MUMPS library \cite{MUMPS,MUMPS1,MUMPS2}.
For simplicity, we neglect the other parts of the computation, such as applications of
the prolongation and restrictions operators or the multiplication of a vector by
the (total) matrix $\mA$ which is executed in iterative solvers.
Moreover, we assume that
we have enough computer cores and that
each algebraic system from \eqref{flops0}
is solved by one core
so that all independent algebraic systems can be solved in parallel. 
We measure computational costs by
\begin{enumerate}[label=({\roman*})]
\item $\Fl$ -- the maximum of the number of {\em floating point operations} per one core,
\item $\comm$ -- the number of {\em communication operations} among the cores.
\end{enumerate}


\subsubsection{Floating point operations}

The solution of each linear algebraic system from \eqref{flops0} by MUMPS has two steps:
\begin{enumerate}[label=({\arabic*})]
\item the {\em factorization} of the system, which is carried out using $\flfac(n)$
  floating point operations, where $n$ denotes the size of the system,
\item the {\em assembling} of the solution using $\flass(n)$ floating point operations.
\end{enumerate}
The factorization is performed only once before the start of the iterative solver, whereas
the assembling is performed at each solver iteration.
Both values $\flfac$ and $\flass$ are provided by MUMPS.


We estimate the maximum number of floating point operations per core of the
iterative solver.
Let $n_i:= \dim\Ship$, $i=1,\dots,N$ denote the dimension of the local spaces \eqref{Ship}
and similarly $n_0:=\dim\SHP$, cf.~\eqref{SHP}.
The factorization {\em} of matrices $\mA_i\in\R^{n_i\times n_i}$ can be carried out independently for
each $i=0,\dots,N$, so the maximum number of floating point operations per core is 
\begin{align}
  \label{flops1}
 \FFfac := \max\nolimits_{i=0,\dots,N} \flfac(n_i).
\end{align}

The {\em assembling} of the solutions of the local systems and the global coarse one
can be executed in parallel only for the additive preconditioner $\mNadd^{-1}$.
For the hybrid preconditioner  $\mNhyb^{-1}$, we solve the local systems in parallel
and (twice) the global system sequentially. Hence, the corresponding
maximum number of floating-point operations per core for the assembling is
\begin{align}
  \label{flops2}
  \FFass :=
  \begin{cases}
    \max_{i=0,\dots,N} \flass(n_i) & \mbox{ for } \mNadd^{-1},\\
    \max_{i=1,\dots,N} \flass(n_i)  + 2 \flass(n_0) & \mbox{ for } \mNhyb^{-1},\\
  \end{cases}
\end{align}
at each solver iteration. Let $\iter$ denote the number of iterations of the iterative solver then,
using \eqref{flops1}--\eqref{flops2}, the maximum of the number of floating point operations
per core is
\begin{align}
  \label{flops3}
  \Fl :=\FFfac + \iter\, \FFass.
\end{align}

\subsection{Communication operations} 

We are aware that the following considerations exhibit significant simplification and that
the number of communication operations strongly depends on the implementation. We assume that
all vectors appearing in the computations are stored in copies at each computer core and that
each matrix $\mA_i$, $i=0,\dots,N$ is allocated only at one processor. However, for the
hybrid operator case, we assume that the ``coarse'' matrix $\mA_0$ is stored in copies at each
core. This causes a (small) increase of the memory requirements, but
the same coarse global problem can be solved at each core, which keeps
the maximum of the floating point operations, whereas reduces the communication among
the cores.

Therefore,
all communication operators among the cores are given by transmission
of the vectors of $n$ (=the size of $\mA$). Using, for example, the tree-based algorithm,
this transmission is proportional to $\log_2{N}$, where $N$ is the number of cores.
Hence, the number of {\em communication operations} can be estimated as
\begin{align}
  \label{comm5}
  \comm =  \iter \,{n}\,\log_2{N},
\end{align}
where $\iter$ is the number of iterations of the iterative solver.
In contrast to floating-point iterations, the number of communications increases with the
number of subdomains (= the number of used cores).

We recall that both \eqref{flops3} and \eqref{comm5} have only an informative character, since
some parts of iterative solvers (multiplication by $\mA$,  application of $\mRi$, $\mRiT$)
are not considered. However, they provide additional information related to computational
costs.

\section{Numerical examples}
\label{sec:numer}

In this section, we present the numerical study of the convergence of the preconditioners
presented in Section~\ref{sec:DDM}. The subdomains $\Om_i$, $i=1,\dots,N$ and 
the coarse meshes $\ThH$ are generated by METIS \cite{metis}.
The aim is to show
\begin{itemize}
\item the {\em weak scalability} of the iterative methods, i.e,
 the computational costs are ideally
 constant for a fixed ratio between the size of the problem and the number of computer cores, 
\item the comparison of efficiency of the additive and hybrid Schwarz preconditioners,
\item the dependence of the convergence on the the diffusivity variation,
\item the dependence of the convergence on
  the domain decomposition respecting and non-respecting the material interfaces.
\end{itemize}
We performed computations using a sequence of (quasi) uniform
meshes with increasing numbers of elements $\#\Th$ and
the number of subdomains $N$ is chosen such that the number of elements
within each $\Om_i$, $i=1,\dots,N$ is (approximately constant). Namely we keep the
ratio $\#\Th/N\approx 100$ and $\#\Th/N\approx \numf{1000}$.
The coarse mesh $\ThH$ is chosen such that each $\cK$ is just one
subdomain $\Om_i$, $i=1,\dots,N$, or each $\Om_i$ is additionally divided into several $\cK$.

Linear systems \eqref{alg1} are solved using the conjugate gradient (CG) method with
the preconditioners $\mNadd^{-1}$ and $\mNhyb^{-1}$ given by \eqref{repre2}
and \eqref{hyb4}, respectively.
The CG algorithm is stopped when the relative preconditioned residual $\rrel$ fulfills
\begin{align}
  \label{alg:stop}
  \rrel^\ell := \| \mN^{-1}(\mA \bku^\ell - \bkg)\|/ \| \mN^{-1}(\mA \bku^0 -\bkg)\| \leq \TOL,
\end{align}
where $\mN^{-1}$ denotes a preconditioner and $\TOL>0$ is the prescribed user tolerance.
We employ $P_p$, $p=1,2,3$ polynomial approximations for each case.

\subsection{Laplace problem}
\label{sec:lapl}
First, we consider a simple toy example 
\begin{align}
  \label{lapl}
  -\Delta u = -2x_1(1-x_1) - 2x_2(1-x_2) \qquad \mbox{in } \Om=(0,1)^2, 
\end{align}
with the homogeneous Dirichlet boundary condition on $\Gamma$, which gives the exact solution
$u=x_1(1-x_1)x_2(1-x_2)$. 
The initial approximation $\bku^0$ corresponds to a highly oscillating function
(namely $u=\sum_{i,j=1}^3 \sin(2\pi i x_1) \sin(2\pi j x_2)$) in order to avoid
a possible superconvergence due to the presence of particular frequency modes, cf.
\cite{GanderTalk}. We set $\TOL=10^{-12}$ in \eqref{alg:stop}.

The results achieved are given in Tables~\ref{tab:Lapl}--\ref{tab:Lapl3}, where we present
\begin{itemize}
  \renewcommand\labelitemi{--}
\item $\#\Th$ -- the number of elements of the fine mesh,
\item $N$ -- the number of subdomains $\Om_i$,
\item  $\#\Thi:=\#\Th/N$ -- the average number of elements in $\Om_i$, $i=1,\dots,N$,
\item $\#\ThH$ -- the number of elements of the coarse  mesh,
\item $\iter$ -- the number of (preconditioned) CG iterations necessary to achieve \eqref{alg:stop},
\item  $\MFl$ -- the number of floating point operations given by \eqref{flops3},
  $\MFl = 10^6\Fl$,
\item  $\Mcomm$ -- the number of communication operations given by  \eqref{comm5}, $\Mcomm = 10^6\comm$.
\end{itemize}
Table~\ref{tab:Lapl} shows the results corresponding to
$\#\Th/N = \#\Thi\approx 100$ and
each subdomain $\Om_i$, $i=1,\dots,N$ is just one element $\cK\in\ThH$.
Furthermore, Table~\ref{tab:Lapl2} contains the results with
$\#\Th/N = \#\Thi\approx 1000$ and each subdomain $\Om_i$ consists of 1, 5, and 10
coarse elements $\cK\in\ThH$. 

\begin{figure}
  \vertical{\footnotesize\hspace{4mm} additive  $\mNadd^{-1}$}
  \includegraphics[width=0.315\textwidth]{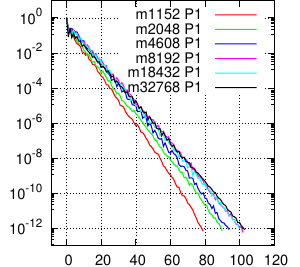}
  \includegraphics[width=0.315\textwidth]{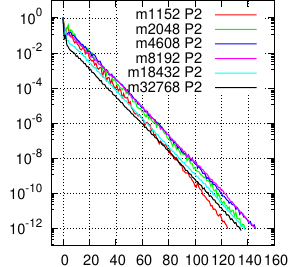}
  \includegraphics[width=0.315\textwidth]{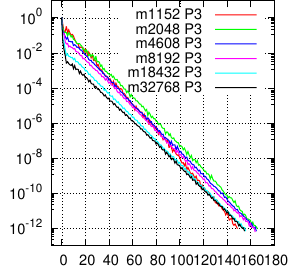}

  \vspace{2mm}
  
  \vertical{\footnotesize\hspace{4mm} hybrid  $\mNhyb^{-1}$}
  \includegraphics[width=0.315\textwidth]{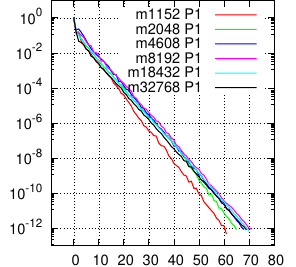}
  \includegraphics[width=0.315\textwidth]{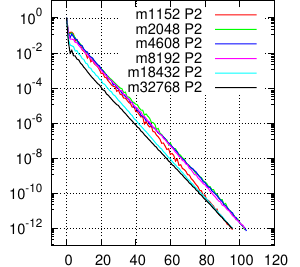}
  \includegraphics[width=0.315\textwidth]{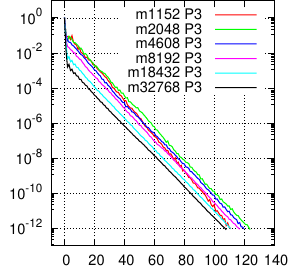}
  \caption{Laplace problem \eqref{lapl}, convergence of CG method
    ($\rrel^\ell$ for $\ell=0,1,\dots$, cf.~\eqref{alg:stop})  for preconditioners
    $\mNadd^{-1}$ and $\mNhyb^{-1}$ with the subdomains splitting $\#\Th/N\approx100$.}
  \label{fig:Lapl}
\end{figure}

Moreover, Table~\ref{tab:Lapl3} presents the results for the finest mesh $\Th$ having
$\numf{32768}$ elements, the number of sub-domains is $N=8,\ 16,\ 32,\ 64$, and
the number of elements of the coarse mesh is fixed to 128
(however, the coarse meshes are not the same due to kind of construction).
Finally, Figure~\ref{fig:Lapl} shows the convergence of the CG method
for the setting from Table~\ref{tab:Lapl}, 
namely the dependence of $\rrel^\ell$ with respect to $\ell=0,1,\dots$, cf.~\eqref{alg:stop}.
We observe the following.
\begin{itemize}
\item Tables~\ref{tab:Lapl}--\ref{tab:Lapl2}:
  For each particular $p$ and particular preconditioner, the number of CG iterations is almost
  constant for increasing $\#\Th$, which implies the weak scalability of both preconditioners
  provided that the computational costs of the coarse solver are neglected.
  This is not the case for the finer meshes in Table~\ref{tab:Lapl} where $\#\ThH > \#\Thi$ and thus the number of $\Fl$ increases with increasing size of the problem.
  However, the results of Table~\ref{tab:Lapl2}, when $\#\ThH \ll \#\Thi$, support the scalability in terms of $\Fl$.
\item Table~\ref{tab:Lapl3}: When $\#\ThH$ and $\#\Thi$ are fixed
(i.e., $h$ and $H$ are fixed, also), we observe a slight increase in the number of 
$\iter$ for an increasing number of subdomains $N$. Consequently, the number of $\Fl$ is reduced by
  factor at least two since the parallelism can be used more effectively. On the other hand,
  the number of communication operations $\comm$ is increasing and prolonging the computation.
  The balance of these two aspects (by the optimal choice of the number of subdomains) is open
  and will be the subject of further research.
\item Tables~\ref{tab:Lapl}--\ref{tab:Lapl3}:
the performance of the hybrid preconditioner $\mNhyb^{-1}$ saves about 25 -- 30\% of the CG iterations compared to the additive preconditioner $\mNadd^{-1}$.
  The value 25 -- 30\% also exhibits a potential benefit in computational time if the computational costs of the coarse solver are negligible (e.g., $\#\ThH  \lesssim \#\Thi$).
\item  Figure~\ref{fig:Lapl}: All graphs showing the convergence of the residual with respect to the
  number of CG iterations support the scalability of the methods. Moreover,
for the $P_2$ and $P_3$ approximation, we observe that the tolerance level $\TOL=10^{-6}$ in \eqref{alg:stop} is achieved using the smaller number of CG iterations for the increasing size of the problem.
\end{itemize}

\begin{table}
  \caption{Laplace problem \eqref{lapl}, convergence of CG method with preconditioners
    $\mNadd^{-1}$ and $\mNhyb^{-1}$, the subdomains splitting $\#\Th/N\approx100$.}
  \label{tab:Lapl}
  \begin{center}
   {\footnotesize
      \tabcolsep 6pt
  \begin{tabular}{cccc|rrr|rrr|rrr}
  \hline
  \multicolumn{4}{c|}{ {\small\bf additive} $\mNadd^{-1}$     }  &  \multicolumn{3}{c|}{$p=1$} & \multicolumn{3}{c|}{$p=2$}&   \multicolumn{3}{c}{$p=3$}\\
  \hline
     $ \# \Th $   &
     $ N $ & $ \# \Thi $   &
     $ \# \ThH $ & 
     $ \iter $ & $ \MFl $ & $ \Mcomm $ & 
     $ \iter $ & $ \MFl $ & $ \Mcomm $ & 
     $ \iter $ & $ \MFl $ & $ \Mcomm $ \\ 
  \hline
    1152 &       11 &      104 &       11
& {\bf       79}  &        0.9 &        0.9
& {\bf      125}  &        5.9 &        3.0
& {\bf      149}  &       21.1 &        5.9
 \\
    2048 &       20 &      102 &       20
& {\bf       90}  &        0.9 &        2.4
& {\bf      138}  &        5.6 &        7.3
& {\bf      166}  &       20.5 &       14.7
 \\
    4608 &       46 &      100 &       46
& {\bf       94}  &        0.9 &        7.2
& {\bf      146}  &        5.8 &       22.3
& {\bf      165}  &       20.7 &       42.0
 \\
    8192 &       81 &      101 &       81
& {\bf      102}  &        1.2 &       15.9
& {\bf      146}  &        6.6 &       45.5
& {\bf      164}  &       25.3 &       85.2
 \\
   18432 &      184 &      100 &      184
& {\bf      100}  &        4.0 &       41.6
& {\bf      138}  &       24.4 &      114.8
& {\bf      156}  &       82.9 &      216.3
 \\
   32768 &      327 &      100 &      327
& {\bf      103}  &        9.5 &       84.6
& {\bf      135}  &       53.8 &      221.7
& {\bf      155}  &      198.3 &      424.3
 \\
  \hline
  \end{tabular}
 }

  \vspace{2mm}
  
   {\footnotesize
      \tabcolsep 6pt
  \begin{tabular}{cccc|rrr|rrr|rrr}
  \hline
  \multicolumn{4}{c|}{{\small\bf hybrid} $\mNhyb^{-1}$        }  &  \multicolumn{3}{c|}{$p=1$} & \multicolumn{3}{c|}{$p=2$}&   \multicolumn{3}{c}{$p=3$}\\
  \hline
     $ \# \Th $   &
     $ N $ & $ \# \Thi $   &
     $ \# \ThH $ & 
     $ \iter $ & $ \MFl $ & $ \Mcomm $ & 
     $ \iter $ & $ \MFl $ & $ \Mcomm $ & 
     $ \iter $ & $ \MFl $ & $ \Mcomm $ \\ 
  \hline
    1152 &       11 &      104 &       11
& {\bf       61}  &        0.7 &        0.7
& {\bf       96}  &        5.1 &        2.3
& {\bf      111}  &       17.9 &        4.4
 \\
    2048 &       20 &      102 &       20
& {\bf       65}  &        0.9 &        1.7
& {\bf      103}  &        5.5 &        5.5
& {\bf      124}  &       20.1 &       11.0
 \\
    4608 &       46 &      100 &       46
& {\bf       68}  &        1.3 &        5.2
& {\bf      104}  &        7.7 &       15.9
& {\bf      120}  &       27.7 &       30.5
 \\
    8192 &       81 &      101 &       81
& {\bf       71}  &        2.3 &       11.1
& {\bf      103}  &       12.4 &       32.1
& {\bf      117}  &       44.5 &       60.8
 \\
   18432 &      184 &      100 &      184
& {\bf       70}  &        5.9 &       29.1
& {\bf       96}  &       34.8 &       79.9
& {\bf      111}  &      115.8 &      153.9
 \\
   32768 &      327 &      100 &      327
& {\bf       69}  &       12.7 &       56.7
& {\bf       96}  &       73.2 &      157.7
& {\bf      108}  &      254.8 &      295.6
 \\
  \hline
  \end{tabular}
 }

  \end{center}
\end{table}

\begin{table}
  \caption{Laplace problem \eqref{lapl}, convergence of CG method with preconditioners
    $\mNadd^{-1}$ and $\mNhyb^{-1}$, the subdomains splitting $\#\Th/N\approx\numf{1000}$.}
  \label{tab:Lapl2}
  \begin{center}
    
     {\footnotesize
      \tabcolsep 6pt
  \begin{tabular}{cccc|rrr|rrr|rrr}
  \hline
  \multicolumn{4}{c|}{ {\small\bf additive} $\mNadd^{-1}$     }  &  \multicolumn{3}{c|}{$p=1$} & \multicolumn{3}{c|}{$p=2$}&   \multicolumn{3}{c}{$p=3$}\\
  \hline
     $ \# \Th $   &
     $ N $ & $ \# \Thi $   &
     $ \# \ThH $ & 
     $ \iter $ & $ \MFl $ & $ \Mcomm $ & 
     $ \iter $ & $ \MFl $ & $ \Mcomm $ & 
     $ \iter $ & $ \MFl $ & $ \Mcomm $ \\ 
  \hline
    8192 &        8 &     1024 &        8
& {\bf      123}  &       33.0 &        9.1
& {\bf      200}  &      249.9 &       29.5
& {\bf      241}  &      883.5 &       59.2
 \\
   18432 &       18 &     1024 &       18
& {\bf      157}  &       46.1 &       36.2
& {\bf      227}  &      270.4 &      104.7
& {\bf      274}  &      999.3 &      210.6
 \\
   32768 &       32 &     1024 &       32
& {\bf      166}  &       45.3 &       81.6
& {\bf      243}  &      301.0 &      238.9
& {\bf      284}  &     1039.8 &      465.3
 \\
  \hline
  \hline
    8192 &        8 &     1024 &       40
& {\bf       89}  &       27.2 &        6.6
& {\bf      138}  &      199.8 &       20.3
& {\bf      166}  &      733.7 &       40.8
 \\
   18432 &       18 &     1024 &       90
& {\bf       99}  &       31.8 &       22.8
& {\bf      143}  &      202.6 &       65.9
& {\bf      167}  &      737.7 &      128.4
 \\
   32768 &       32 &     1024 &      160
& {\bf      103}  &       30.5 &       50.6
& {\bf      138}  &      209.3 &      135.7
& {\bf      168}  &      753.1 &      275.2
 \\
  \hline
  \hline
    8192 &        8 &     1024 &       80
& {\bf       78}  &       24.5 &        5.8
& {\bf      117}  &      170.6 &       17.3
& {\bf      137}  &      642.8 &       33.7
 \\
   18432 &       18 &     1024 &      180
& {\bf       82}  &       29.2 &       18.9
& {\bf      112}  &      173.2 &       51.7
& {\bf      134}  &      665.4 &      103.0
 \\
   32768 &       32 &     1024 &      320
& {\bf       86}  &       28.8 &       42.3
& {\bf      114}  &      175.5 &      112.1
& {\bf      136}  &      670.6 &      222.8
 \\
  \hline
  \end{tabular}
 }

    \vspace{2mm}
    
     {\footnotesize
      \tabcolsep 6pt
  \begin{tabular}{cccc|rrr|rrr|rrr}
  \hline
  \multicolumn{4}{c|}{{\small\bf hybrid} $\mNhyb^{-1}$        }  &  \multicolumn{3}{c|}{$p=1$} & \multicolumn{3}{c|}{$p=2$}&   \multicolumn{3}{c}{$p=3$}\\
  \hline
     $ \# \Th $   &
     $ N $ & $ \# \Thi $   &
     $ \# \ThH $ & 
     $ \iter $ & $ \MFl $ & $ \Mcomm $ & 
     $ \iter $ & $ \MFl $ & $ \Mcomm $ & 
     $ \iter $ & $ \MFl $ & $ \Mcomm $ \\ 
  \hline
    8192 &        8 &     1024 &        8
& {\bf      100}  &       28.3 &        7.4
& {\bf      152}  &      208.0 &       22.4
& {\bf      185}  &      751.5 &       45.5
 \\
   18432 &       18 &     1024 &       18
& {\bf      117}  &       36.9 &       27.0
& {\bf      173}  &      224.6 &       79.8
& {\bf      205}  &      834.0 &      157.6
 \\
   32768 &       32 &     1024 &       32
& {\bf      126}  &       37.1 &       61.9
& {\bf      181}  &      247.3 &      177.9
& {\bf      216}  &      882.4 &      353.9
 \\
  \hline
  \hline
    8192 &        8 &     1024 &       40
& {\bf       65}  &       22.6 &        4.8
& {\bf      103}  &      170.3 &       15.2
& {\bf      125}  &      644.0 &       30.7
 \\
   18432 &       18 &     1024 &       90
& {\bf       72}  &       27.4 &       16.6
& {\bf      104}  &      176.8 &       48.0
& {\bf      123}  &      658.9 &       94.5
 \\
   32768 &       32 &     1024 &      160
& {\bf       74}  &       28.3 &       36.4
& {\bf      102}  &      195.8 &      100.3
& {\bf      121}  &      703.4 &      198.2
 \\
  \hline
  \hline
    8192 &        8 &     1024 &       80
& {\bf       55}  &       20.9 &        4.1
& {\bf       83}  &      146.6 &       12.2
& {\bf      101}  &      578.3 &       24.8
 \\
   18432 &       18 &     1024 &      180
& {\bf       58}  &       27.7 &       13.4
& {\bf       81}  &      163.6 &       37.4
& {\bf       97}  &      638.0 &       74.6
 \\
   32768 &       32 &     1024 &      320
& {\bf       59}  &       30.5 &       29.0
& {\bf       82}  &      190.9 &       80.6
& {\bf       96}  &      713.1 &      157.3
 \\
  \hline
  \end{tabular}
 }

  \end{center}
\end{table}

\begin{table}
  \caption{Laplace problem \eqref{lapl}, convergence of CG method
    with preconditioners
    $\mNadd^{-1}$ and $\mNhyb^{-1}$ for the finest mesh with fixed $\#\ThH$
    using various domain decomposition.}
  \label{tab:Lapl3}
  \begin{center}
     {\footnotesize
      \tabcolsep 6pt
  \begin{tabular}{cccc|rrr|rrr|rrr}
  \hline
  \multicolumn{4}{c|}{ {\small\bf additive} $\mNadd^{-1}$     }  &  \multicolumn{3}{c|}{$p=1$} & \multicolumn{3}{c|}{$p=2$}&   \multicolumn{3}{c}{$p=3$}\\
  \hline
     $ \# \Th $   &
     $ N $ & $ \# \Thi $   &
     $ \# \ThH $ & 
     $ \iter $ & $ \MFl $ & $ \Mcomm $ & 
     $ \iter $ & $ \MFl $ & $ \Mcomm $ & 
     $ \iter $ & $ \MFl $ & $ \Mcomm $ \\ 
  \hline
   32768 &        8 &     4096 &      128
& {\bf       94}  &      198.6 &       27.7
& {\bf      144}  &     1370.6 &       84.9
& {\bf      166}  &     5338.8 &      163.2
 \\
   32768 &       16 &     2048 &      128
& {\bf       99}  &       80.1 &       38.9
& {\bf      151}  &      549.4 &      118.8
& {\bf      174}  &     2135.1 &      228.1
 \\
   32768 &       32 &     1024 &      128
& {\bf      107}  &       31.2 &       52.6
& {\bf      151}  &      215.6 &      148.4
& {\bf      176}  &      770.6 &      288.4
 \\
   32768 &       64 &      512 &      128
& {\bf      114}  &       13.2 &       67.2
& {\bf      159}  &       74.5 &      187.6
& {\bf      182}  &      294.0 &      357.8
 \\
  \hline
  \end{tabular}
 }

    \vspace{2mm}
    
     {\footnotesize
      \tabcolsep 6pt
  \begin{tabular}{cccc|rrr|rrr|rrr}
  \hline
  \multicolumn{4}{c|}{{\small\bf hybrid} $\mNhyb^{-1}$        }  &  \multicolumn{3}{c|}{$p=1$} & \multicolumn{3}{c|}{$p=2$}&   \multicolumn{3}{c}{$p=3$}\\
  \hline
     $ \# \Th $   &
     $ N $ & $ \# \Thi $   &
     $ \# \ThH $ & 
     $ \iter $ & $ \MFl $ & $ \Mcomm $ & 
     $ \iter $ & $ \MFl $ & $ \Mcomm $ & 
     $ \iter $ & $ \MFl $ & $ \Mcomm $ \\ 
  \hline
   32768 &        8 &     4096 &      128
& {\bf       69}  &      169.3 &       20.3
& {\bf      104}  &     1181.7 &       61.3
& {\bf      123}  &     4784.0 &      120.9
 \\
   32768 &       16 &     2048 &      128
& {\bf       75}  &       69.7 &       29.5
& {\bf      112}  &      479.8 &       88.1
& {\bf      130}  &     1920.0 &      170.4
 \\
   32768 &       32 &     1024 &      128
& {\bf       80}  &       28.4 &       39.3
& {\bf      111}  &      194.3 &      109.1
& {\bf      131}  &      716.5 &      214.6
 \\
   32768 &       64 &      512 &      128
& {\bf       81}  &       12.9 &       47.8
& {\bf      115}  &       76.5 &      135.7
& {\bf      133}  &      296.8 &      261.5
 \\
  \hline
  \end{tabular}
 }

  \end{center}
\end{table}

\subsection{Alternator (linearized)}
\label{sec:alter}
\begin{figure} [t]
  \begin{center}
%
    \includegraphics[width=0.31\textwidth]{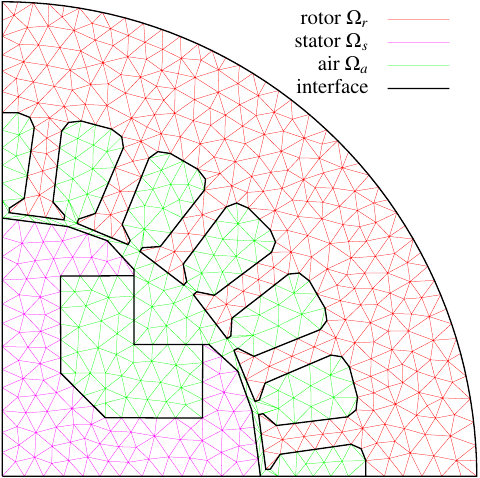}
    \hspace{0.01\textwidth}
    \includegraphics[width=0.31\textwidth]{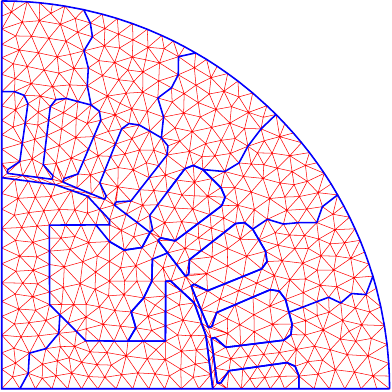}
    \hspace{0.01\textwidth}
    \includegraphics[width=0.31\textwidth]{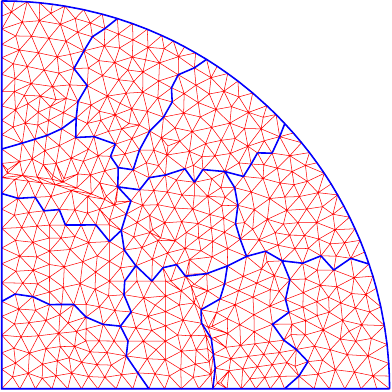}
\end{center}
  \caption{Alternator, the computational domain $\Om$ with its components (left)
    and the computational meshes $\Th$ (red) and $\ThH$ (blue) with $\#\Th=1112$ \& $N=\#\ThH=11$:
    respecting (center) and non-respecting (right) the component interfaces.}
  \label{fig:alter_geom}
\end{figure}

This example exhibits a linearized variant of
the magnetic state in the cross section of an alternator from Section~\ref{sec:alterN},
which originates from \cite{Glowinski74}.
The 
computational domain $\Om$ (one quarter of the circle alternator)
consists of several components represented by
the stator ($\Omega_s$) and the rotor ($\Omega_r$)
with a gap filled by air ($\Omega_a$);
see Figure~\ref{fig:alter_geom}, left, where the geometry of the domain is shown.
We consider problem \eqref{prob1a} with 
$f =  5\cdot 10^{4}$ 
and $\K = \K(x) $ such that
\begin{align}
  \label{alter1}
  \K(x) = 
  \begin{cases}
    \mu_0^{-1} & \mbox{ for } x\in \Om_a, \\
    \mu_1^{-1} &  \mbox{ for } x\in \Om_s\cup\Om_r, \\
  \end{cases}
\end{align}
where $\mu_0 =  1.256\cdot 10^{-6}$, $\mu_1 = \zeta \mu_0$, $\zeta=\barK/\ubarK > 1$. Particularly,
we use the values $\zeta=100$ and $\zeta=\numf{10000}$.
We prescribe 
$\nabla u\cdot\bkn = 0$ on $\gomN:=(0,1)\times\{0\}\cup \{0\} \times (0,1)$ and
$u=0$  
on the rest of the boundary. 

In accordance with Remark~\ref{rem:materialDD}, we consider two types of domain decomposition.
\begin{itemize}
\item {\em respecting} the component interfaces:
each subdomain $\Om_i$, $i=1,\dots, N$ is subset either of $\Om_s$ or $\Om_r$ or $\Om_a$,
i.e., $\Om_i$ does not contain any component interface, 
\item {\em non-respecting} the component interfaces:
  the domain partition does not fulfill the previous property and the domain decomposition
  is carried out without any restriction related to component interfaces.
\end{itemize}
Figure~\ref{fig:alter_geom}, center and right, shows an example.

In the same way as in Section~\ref{sec:lapl},
we performed computations using a sequence of (quasi) uniform meshes
and domain decomposition
such that the number of elements within each subdomain is (approximately) fixed.
We present the quantities $\#\Th$, $\#\Thi$, $N$, and $\#\ThH$ together with
the number of (preconditioned) CG iterations $\iter$ necessary to achieve \eqref{alg:stop}
(with $\TOL=10^{-10}$),
and the computational costs in $\MFl$ and $\Mcomm$.
Tables~\ref{tab:alter} and \ref{tab:alter2} show the results for
$\zeta = \mu_1/\mu_0=100$ in \eqref{alter1} for $\#\Th/N\approx\numf{100}$ and
$\#\Th/N\approx\numf{1000}$, respectively.
Moreover, Tables~\ref{tab:alter4} and \ref{tab:alter5} present results
corresponding to $\zeta=\numf{10000}$ for $\#\Th/N\approx\numf{100}$ and
$\#\Th/N\approx\numf{1000}$, respectively.

\begin{table}
  \caption{Linearized alternator \eqref{alter1} with $\zeta=100$,
    convergence of CG method with preconditioners
    $\mNadd^{-1}$ and $\mNhyb^{-1}$ with the subdomains splitting $\#\Th/N\approx100$.}
  \label{tab:alter}
  \begin{center}

    {
      domain decomposition {\bf non-respecting} component interfaces}
    
     {\footnotesize
      \tabcolsep 6pt
  \begin{tabular}{cccc|rrr|rrr|rrr}
  \hline
  \multicolumn{4}{c|}{ {\small\bf additive} $\mNadd^{-1}$     }  &  \multicolumn{3}{c|}{$p=1$} & \multicolumn{3}{c|}{$p=2$}&   \multicolumn{3}{c}{$p=3$}\\
  \hline
     $ \# \Th $   &
     $ N $ & $ \# \Thi $   &
     $ \# \ThH $ & 
     $ \iter $ & $ \MFl $ & $ \Mcomm $ & 
     $ \iter $ & $ \MFl $ & $ \Mcomm $ & 
     $ \iter $ & $ \MFl $ & $ \Mcomm $ \\ 
  \hline
    1112 &       11 &      101 &       11
& {\bf       87}  &        0.9 &        1.0
& {\bf      149}  &        5.8 &        3.4
& {\bf      157}  &       18.4 &        6.0
 \\
    2025 &       20 &      101 &       20
& {\bf      112}  &        1.1 &        2.9
& {\bf      133}  &        5.6 &        7.0
& {\bf      168}  &       21.5 &       14.7
 \\
    3832 &       38 &      100 &       38
& {\bf      110}  &        1.1 &        6.6
& {\bf      143}  &        6.2 &       17.3
& {\bf      172}  &       21.5 &       34.6
 \\
    7549 &       75 &      100 &       75
& {\bf      135}  &        1.5 &       19.0
& {\bf      163}  &        6.9 &       46.0
& {\bf      186}  &       25.3 &       87.5
 \\
   14937 &      149 &      100 &      149
& {\bf      121}  &        3.9 &       39.1
& {\bf      155}  &       19.2 &      100.3
& {\bf      162}  &       65.0 &      174.7
 \\
   29819 &      298 &      100 &      298
& {\bf      112}  &        9.4 &       82.3
& {\bf      155}  &       54.7 &      227.9
& {\bf      159}  &      193.8 &      389.7
 \\
  \hline
  \end{tabular}
 }

     {\footnotesize
      \tabcolsep 6pt
  \begin{tabular}{cccc|rrr|rrr|rrr}
  \hline
  \multicolumn{4}{c|}{{\small\bf hybrid} $\mNhyb^{-1}$        }  &  \multicolumn{3}{c|}{$p=1$} & \multicolumn{3}{c|}{$p=2$}&   \multicolumn{3}{c}{$p=3$}\\
  \hline
     $ \# \Th $   &
     $ N $ & $ \# \Thi $   &
     $ \# \ThH $ & 
     $ \iter $ & $ \MFl $ & $ \Mcomm $ & 
     $ \iter $ & $ \MFl $ & $ \Mcomm $ & 
     $ \iter $ & $ \MFl $ & $ \Mcomm $ \\ 
  \hline
    1112 &       11 &      101 &       11
& {\bf       62}  &        0.7 &        0.7
& {\bf      112}  &        5.0 &        2.6
& {\bf      107}  &       14.9 &        4.1
 \\
    2025 &       20 &      101 &       20
& {\bf       77}  &        1.0 &        2.0
& {\bf       95}  &        5.4 &        5.0
& {\bf      106}  &       18.2 &        9.3
 \\
    3832 &       38 &      100 &       38
& {\bf       76}  &        1.3 &        4.6
& {\bf       97}  &        6.8 &       11.7
& {\bf      106}  &       23.7 &       21.3
 \\
    7549 &       75 &      100 &       75
& {\bf       95}  &        2.9 &       13.4
& {\bf      110}  &       12.8 &       31.0
& {\bf      119}  &       42.3 &       56.0
 \\
   14937 &      149 &      100 &      149
& {\bf       85}  &        6.0 &       27.5
& {\bf       97}  &       26.5 &       62.8
& {\bf       99}  &       85.2 &      106.8
 \\
   29819 &      298 &      100 &      298
& {\bf       68}  &       11.7 &       50.0
& {\bf       97}  &       68.2 &      142.6
& {\bf       94}  &      225.0 &      230.4
 \\
  \hline
  \end{tabular}
 }

    \vspace{2mm}
    
    {
       domain decomposition {\bf respecting} component interfaces}
    
     {\footnotesize
      \tabcolsep 6pt
  \begin{tabular}{cccc|rrr|rrr|rrr}
  \hline
  \multicolumn{4}{c|}{ {\small\bf additive} $\mNadd^{-1}$     }  &  \multicolumn{3}{c|}{$p=1$} & \multicolumn{3}{c|}{$p=2$}&   \multicolumn{3}{c}{$p=3$}\\
  \hline
     $ \# \Th $   &
     $ N $ & $ \# \Thi $   &
     $ \# \ThH $ & 
     $ \iter $ & $ \MFl $ & $ \Mcomm $ & 
     $ \iter $ & $ \MFl $ & $ \Mcomm $ & 
     $ \iter $ & $ \MFl $ & $ \Mcomm $ \\ 
  \hline
    1112 &       11 &      101 &       11
& {\bf      350}  &        2.7 &        4.0
& {\bf      327}  &       10.4 &        7.5
& {\bf      474}  &       42.5 &       18.2
 \\
    2025 &       19 &      106 &       19
& {\bf      212}  &        2.4 &        5.5
& {\bf      287}  &       13.3 &       14.8
& {\bf      318}  &       42.1 &       27.4
 \\
    3832 &       38 &      100 &       38
& {\bf      196}  &        1.8 &       11.8
& {\bf      309}  &       11.7 &       37.3
& {\bf      323}  &       36.0 &       65.0
 \\
    7549 &       75 &      100 &       75
& {\bf      204}  &        2.0 &       28.8
& {\bf      292}  &       11.6 &       82.4
& {\bf      309}  &       38.8 &      145.3
 \\
   14937 &      149 &      100 &      149
& {\bf      191}  &        4.7 &       61.8
& {\bf      237}  &       23.8 &      153.3
& {\bf      242}  &       75.3 &      261.0
 \\
   29819 &      298 &      100 &      298
& {\bf      176}  &       11.8 &      129.4
& {\bf      181}  &       54.3 &      266.2
& {\bf      196}  &      193.4 &      480.4
 \\
  \hline
  \end{tabular}
 }

     {\footnotesize
      \tabcolsep 6pt
  \begin{tabular}{cccc|rrr|rrr|rrr}
  \hline
  \multicolumn{4}{c|}{{\small\bf hybrid} $\mNhyb^{-1}$        }  &  \multicolumn{3}{c|}{$p=1$} & \multicolumn{3}{c|}{$p=2$}&   \multicolumn{3}{c}{$p=3$}\\
  \hline
     $ \# \Th $   &
     $ N $ & $ \# \Thi $   &
     $ \# \ThH $ & 
     $ \iter $ & $ \MFl $ & $ \Mcomm $ & 
     $ \iter $ & $ \MFl $ & $ \Mcomm $ & 
     $ \iter $ & $ \MFl $ & $ \Mcomm $ \\ 
  \hline
    1112 &       11 &      101 &       11
& {\bf      267}  &        2.3 &        3.1
& {\bf      227}  &        8.1 &        5.2
& {\bf      315}  &       32.0 &       12.1
 \\
    2025 &       19 &      106 &       19
& {\bf      145}  &        1.9 &        3.7
& {\bf      203}  &       11.1 &       10.5
& {\bf      202}  &       32.7 &       17.4
 \\
    3832 &       38 &      100 &       38
& {\bf      141}  &        2.2 &        8.5
& {\bf      199}  &       12.8 &       24.0
& {\bf      196}  &       36.9 &       39.4
 \\
    7549 &       75 &      100 &       75
& {\bf      133}  &        3.4 &       18.8
& {\bf      189}  &       18.3 &       53.3
& {\bf      192}  &       56.3 &       90.3
 \\
   14937 &      149 &      100 &      149
& {\bf      126}  &        7.2 &       40.8
& {\bf      153}  &       35.2 &       99.0
& {\bf      145}  &      101.5 &      156.4
 \\
   29819 &      298 &      100 &      298
& {\bf      113}  &       15.8 &       83.1
& {\bf      107}  &       66.0 &      157.4
& {\bf      118}  &      232.3 &      289.2
 \\
  \hline
  \end{tabular}
 }

  \end{center}
\end{table}

\begin{table}
  \caption{Linearized alternator \eqref{alter1} with $\zeta=100$,
    convergence of CG method with preconditioners
    $\mNadd^{-1}$ and $\mNhyb^{-1}$ with the subdomains splitting $\#\Th/N\approx\numf{1000}$.}
  \label{tab:alter2}
  \begin{center}

    {
       domain decomposition {\bf non-respecting} component interfaces}
    
     {\footnotesize
      \tabcolsep 6pt
  \begin{tabular}{cccc|rrr|rrr|rrr}
  \hline
  \multicolumn{4}{c|}{ {\small\bf additive} $\mNadd^{-1}$     }  &  \multicolumn{3}{c|}{$p=1$} & \multicolumn{3}{c|}{$p=2$}&   \multicolumn{3}{c}{$p=3$}\\
  \hline
     $ \# \Th $   &
     $ N $ & $ \# \Thi $   &
     $ \# \ThH $ & 
     $ \iter $ & $ \MFl $ & $ \Mcomm $ & 
     $ \iter $ & $ \MFl $ & $ \Mcomm $ & 
     $ \iter $ & $ \MFl $ & $ \Mcomm $ \\ 
  \hline
    7549 &        7 &     1078 &        7
& {\bf      158}  &       36.8 &       10.0
& {\bf      238}  &      250.9 &       30.3
& {\bf      244}  &      898.2 &       51.7
 \\
   14937 &       14 &     1066 &       14
& {\bf      167}  &       41.1 &       28.5
& {\bf      208}  &      243.1 &       71.0
& {\bf      237}  &      812.9 &      134.8
 \\
   29819 &       29 &     1028 &       29
& {\bf      165}  &       38.9 &       71.7
& {\bf      214}  &      238.2 &      186.0
& {\bf      231}  &      843.4 &      334.6
 \\
  \hline
  \end{tabular}
 }

     {\footnotesize
      \tabcolsep 6pt
  \begin{tabular}{cccc|rrr|rrr|rrr}
  \hline
  \multicolumn{4}{c|}{{\small\bf hybrid} $\mNhyb^{-1}$        }  &  \multicolumn{3}{c|}{$p=1$} & \multicolumn{3}{c|}{$p=2$}&   \multicolumn{3}{c}{$p=3$}\\
  \hline
     $ \# \Th $   &
     $ N $ & $ \# \Thi $   &
     $ \# \ThH $ & 
     $ \iter $ & $ \MFl $ & $ \Mcomm $ & 
     $ \iter $ & $ \MFl $ & $ \Mcomm $ & 
     $ \iter $ & $ \MFl $ & $ \Mcomm $ \\ 
  \hline
    7549 &        7 &     1078 &        7
& {\bf      128}  &       30.9 &        8.1
& {\bf      181}  &      205.0 &       23.0
& {\bf      166}  &      704.5 &       35.2
 \\
   14937 &       14 &     1066 &       14
& {\bf      129}  &       33.4 &       22.0
& {\bf      150}  &      193.2 &       51.2
& {\bf      161}  &      638.2 &       91.6
 \\
   29819 &       29 &     1028 &       29
& {\bf      122}  &       31.0 &       53.0
& {\bf      152}  &      189.0 &      132.1
& {\bf      151}  &      659.3 &      218.7
 \\
  \hline
  \end{tabular}
 }

    \vspace{2mm}
    
    {
       domain decomposition {\bf respecting} component interfaces}
    
     {\footnotesize
      \tabcolsep 6pt
  \begin{tabular}{cccc|rrr|rrr|rrr}
  \hline
  \multicolumn{4}{c|}{ {\small\bf additive} $\mNadd^{-1}$     }  &  \multicolumn{3}{c|}{$p=1$} & \multicolumn{3}{c|}{$p=2$}&   \multicolumn{3}{c}{$p=3$}\\
  \hline
     $ \# \Th $   &
     $ N $ & $ \# \Thi $   &
     $ \# \ThH $ & 
     $ \iter $ & $ \MFl $ & $ \Mcomm $ & 
     $ \iter $ & $ \MFl $ & $ \Mcomm $ & 
     $ \iter $ & $ \MFl $ & $ \Mcomm $ \\ 
  \hline
    7549 &        7 &     1078 &        7
& {\bf      337}  &       67.6 &       21.4
& {\bf      524}  &      491.4 &       66.6
& {\bf      606}  &     1657.7 &      128.4
 \\
   14937 &       14 &     1066 &       14
& {\bf      329}  &       68.9 &       56.1
& {\bf      501}  &      493.0 &      170.9
& {\bf      511}  &     1429.6 &      290.6
 \\
   29819 &       29 &     1028 &       29
& {\bf      324}  &       74.8 &      140.8
& {\bf      471}  &      470.3 &      409.4
& {\bf      485}  &     1461.5 &      702.6
 \\
  \hline
  \end{tabular}
 }

     {\footnotesize
      \tabcolsep 6pt
  \begin{tabular}{cccc|rrr|rrr|rrr}
  \hline
  \multicolumn{4}{c|}{{\small\bf hybrid} $\mNhyb^{-1}$        }  &  \multicolumn{3}{c|}{$p=1$} & \multicolumn{3}{c|}{$p=2$}&   \multicolumn{3}{c}{$p=3$}\\
  \hline
     $ \# \Th $   &
     $ N $ & $ \# \Thi $   &
     $ \# \ThH $ & 
     $ \iter $ & $ \MFl $ & $ \Mcomm $ & 
     $ \iter $ & $ \MFl $ & $ \Mcomm $ & 
     $ \iter $ & $ \MFl $ & $ \Mcomm $ \\ 
  \hline
    7549 &        7 &     1078 &        7
& {\bf      244}  &       50.0 &       15.5
& {\bf      385}  &      373.7 &       49.0
& {\bf      442}  &     1270.7 &       93.7
 \\
   14937 &       14 &     1066 &       14
& {\bf      242}  &       52.4 &       41.3
& {\bf      336}  &      352.6 &      114.7
& {\bf      313}  &      979.6 &      178.0
 \\
   29819 &       29 &     1028 &       29
& {\bf      213}  &       52.4 &       92.6
& {\bf      306}  &      330.8 &      266.0
& {\bf      320}  &     1072.9 &      463.5
 \\
  \hline
  \end{tabular}
 }

  \end{center}
\end{table}

\begin{table}
  \caption{Linearized alternator \eqref{alter1} with $\zeta=\numf{10000}$,
    convergence of CG method with preconditioners
    $\mNadd^{-1}$ and $\mNhyb^{-1}$ with the subdomains splitting $\#\Th/N\approx100$.}
  \label{tab:alter4}

  \begin{center}

    {
       domain decomposition {\bf non-respecting} component interfaces}
    
     {\footnotesize
      \tabcolsep 6pt
  \begin{tabular}{cccc|rrr|rrr|rrr}
  \hline
  \multicolumn{4}{c|}{ {\small\bf additive} $\mNadd^{-1}$     }  &  \multicolumn{3}{c|}{$p=1$} & \multicolumn{3}{c|}{$p=2$}&   \multicolumn{3}{c}{$p=3$}\\
  \hline
     $ \# \Th $   &
     $ N $ & $ \# \Thi $   &
     $ \# \ThH $ & 
     $ \iter $ & $ \MFl $ & $ \Mcomm $ & 
     $ \iter $ & $ \MFl $ & $ \Mcomm $ & 
     $ \iter $ & $ \MFl $ & $ \Mcomm $ \\ 
  \hline
    1112 &       11 &      101 &       11
& {\bf      107}  &        1.0 &        1.2
& {\bf      191}  &        7.2 &        4.4
& {\bf      208}  &       23.1 &        8.0
 \\
    2025 &       20 &      101 &       20
& {\bf      140}  &        1.3 &        3.7
& {\bf      178}  &        7.2 &        9.3
& {\bf      238}  &       28.4 &       20.8
 \\
    3832 &       38 &      100 &       38
& {\bf      139}  &        1.3 &        8.4
& {\bf      201}  &        8.4 &       24.3
& {\bf      256}  &       29.8 &       51.5
 \\
    7549 &       75 &      100 &       75
& {\bf      173}  &        1.8 &       24.4
& {\bf      232}  &        9.4 &       65.5
& {\bf      275}  &       34.3 &      129.3
 \\
   14937 &      149 &      100 &      149
& {\bf      163}  &        5.0 &       52.7
& {\bf      234}  &       26.5 &      151.4
& {\bf      271}  &       94.4 &      292.2
 \\
   29819 &      298 &      100 &      298
& {\bf      163}  &       12.9 &      119.8
& {\bf      247}  &       78.3 &      363.2
& {\bf      280}  &      284.6 &      686.2
 \\
  \hline
  \end{tabular}
 }

     {\footnotesize
      \tabcolsep 6pt
  \begin{tabular}{cccc|rrr|rrr|rrr}
  \hline
  \multicolumn{4}{c|}{{\small\bf hybrid} $\mNhyb^{-1}$        }  &  \multicolumn{3}{c|}{$p=1$} & \multicolumn{3}{c|}{$p=2$}&   \multicolumn{3}{c}{$p=3$}\\
  \hline
     $ \# \Th $   &
     $ N $ & $ \# \Thi $   &
     $ \# \ThH $ & 
     $ \iter $ & $ \MFl $ & $ \Mcomm $ & 
     $ \iter $ & $ \MFl $ & $ \Mcomm $ & 
     $ \iter $ & $ \MFl $ & $ \Mcomm $ \\ 
  \hline
    1112 &       11 &      101 &       11
& {\bf       79}  &        0.9 &        0.9
& {\bf      143}  &        6.1 &        3.3
& {\bf      147}  &       19.0 &        5.7
 \\
    2025 &       20 &      101 &       20
& {\bf      102}  &        1.3 &        2.7
& {\bf      126}  &        6.9 &        6.6
& {\bf      159}  &       24.8 &       13.9
 \\
    3832 &       38 &      100 &       38
& {\bf      100}  &        1.6 &        6.0
& {\bf      144}  &        9.7 &       17.4
& {\bf      167}  &       34.7 &       33.6
 \\
    7549 &       75 &      100 &       75
& {\bf      127}  &        3.8 &       17.9
& {\bf      164}  &       18.5 &       46.3
& {\bf      176}  &       59.5 &       82.8
 \\
   14937 &      149 &      100 &      149
& {\bf      121}  &        8.2 &       39.1
& {\bf      158}  &       40.1 &      102.2
& {\bf      158}  &      123.2 &      170.4
 \\
   29819 &      298 &      100 &      298
& {\bf      104}  &       16.9 &       76.5
& {\bf      160}  &      102.8 &      235.3
& {\bf      157}  &      326.0 &      384.8
 \\
  \hline
  \end{tabular}
 }

    \vspace{2mm}
    
    {
       domain decomposition {\bf respecting} component interfaces}
    
     {\footnotesize
      \tabcolsep 6pt
  \begin{tabular}{cccc|rrr|rrr|rrr}
  \hline
  \multicolumn{4}{c|}{ {\small\bf additive} $\mNadd^{-1}$     }  &  \multicolumn{3}{c|}{$p=1$} & \multicolumn{3}{c|}{$p=2$}&   \multicolumn{3}{c}{$p=3$}\\
  \hline
     $ \# \Th $   &
     $ N $ & $ \# \Thi $   &
     $ \# \ThH $ & 
     $ \iter $ & $ \MFl $ & $ \Mcomm $ & 
     $ \iter $ & $ \MFl $ & $ \Mcomm $ & 
     $ \iter $ & $ \MFl $ & $ \Mcomm $ \\ 
  \hline
    1112 &       11 &      101 &       11
& {\bf      412}  &        3.2 &        4.8
& {\bf      414}  &       13.0 &        9.6
& {\bf      518}  &       46.1 &       19.9
 \\
    2025 &       19 &      106 &       19
& {\bf      260}  &        2.9 &        6.7
& {\bf      379}  &       17.3 &       19.6
& {\bf      489}  &       62.2 &       42.1
 \\
    3832 &       38 &      100 &       38
& {\bf      267}  &        2.4 &       16.1
& {\bf      402}  &       15.0 &       48.5
& {\bf      452}  &       48.5 &       90.9
 \\
    7549 &       75 &      100 &       75
& {\bf      265}  &        2.6 &       37.4
& {\bf      366}  &       14.2 &      103.3
& {\bf      395}  &       48.1 &      185.7
 \\
   14937 &      149 &      100 &      149
& {\bf      241}  &        5.8 &       78.0
& {\bf      327}  &       31.4 &      211.6
& {\bf      319}  &       94.0 &      344.0
 \\
   29819 &      298 &      100 &      298
& {\bf      231}  &       15.0 &      169.9
& {\bf      251}  &       70.8 &      369.1
& {\bf      249}  &      229.4 &      610.3
 \\
  \hline
  \end{tabular}
 }

     {\footnotesize
      \tabcolsep 6pt
  \begin{tabular}{cccc|rrr|rrr|rrr}
  \hline
  \multicolumn{4}{c|}{{\small\bf hybrid} $\mNhyb^{-1}$        }  &  \multicolumn{3}{c|}{$p=1$} & \multicolumn{3}{c|}{$p=2$}&   \multicolumn{3}{c}{$p=3$}\\
  \hline
     $ \# \Th $   &
     $ N $ & $ \# \Thi $   &
     $ \# \ThH $ & 
     $ \iter $ & $ \MFl $ & $ \Mcomm $ & 
     $ \iter $ & $ \MFl $ & $ \Mcomm $ & 
     $ \iter $ & $ \MFl $ & $ \Mcomm $ \\ 
  \hline
    1112 &       11 &      101 &       11
& {\bf      350}  &        3.0 &        4.0
& {\bf      301}  &       10.6 &        6.9
& {\bf      348}  &       35.0 &       13.4
 \\
    2025 &       19 &      106 &       19
& {\bf      198}  &        2.6 &        5.1
& {\bf      285}  &       15.2 &       14.7
& {\bf      332}  &       50.7 &       28.6
 \\
    3832 &       38 &      100 &       38
& {\bf      198}  &        3.1 &       11.9
& {\bf      272}  &       17.2 &       32.8
& {\bf      296}  &       53.3 &       59.5
 \\
    7549 &       75 &      100 &       75
& {\bf      184}  &        4.7 &       26.0
& {\bf      247}  &       23.6 &       69.7
& {\bf      255}  &       73.0 &      119.9
 \\
   14937 &      149 &      100 &      149
& {\bf      167}  &        9.5 &       54.0
& {\bf      218}  &       48.6 &      141.1
& {\bf      194}  &      130.2 &      209.2
 \\
   29819 &      298 &      100 &      298
& {\bf      159}  &       21.5 &      116.9
& {\bf      157}  &       91.4 &      230.9
& {\bf      147}  &      274.8 &      360.3
 \\
  \hline
  \end{tabular}
 }

  \end{center}
\end{table}

\begin{table}
  \caption{Linearized alternator \eqref{alter1} with $\zeta=\numf{10000}$,
    convergence of CG method with preconditioners
    $\mNadd^{-1}$ and $\mNhyb^{-1}$ with the subdomains splitting $\#\Th/N\approx\numf{1000}$.}
  \label{tab:alter5}
  \begin{center}

    {
       domain decomposition {\bf non-respecting} component interfaces}
    
     {\footnotesize
      \tabcolsep 6pt
  \begin{tabular}{cccc|rrr|rrr|rrr}
  \hline
  \multicolumn{4}{c|}{ {\small\bf additive} $\mNadd^{-1}$     }  &  \multicolumn{3}{c|}{$p=1$} & \multicolumn{3}{c|}{$p=2$}&   \multicolumn{3}{c}{$p=3$}\\
  \hline
     $ \# \Th $   &
     $ N $ & $ \# \Thi $   &
     $ \# \ThH $ & 
     $ \iter $ & $ \MFl $ & $ \Mcomm $ & 
     $ \iter $ & $ \MFl $ & $ \Mcomm $ & 
     $ \iter $ & $ \MFl $ & $ \Mcomm $ \\ 
  \hline
    7549 &        7 &     1078 &        7
& {\bf      192}  &       43.5 &       12.2
& {\bf      294}  &      296.3 &       37.4
& {\bf      306}  &     1052.5 &       64.8
 \\
   14937 &       14 &     1066 &       14
& {\bf      208}  &       49.6 &       35.5
& {\bf      283}  &      308.3 &       96.6
& {\bf      326}  &     1019.3 &      185.4
 \\
   29819 &       29 &     1028 &       29
& {\bf      211}  &       48.1 &       91.7
& {\bf      306}  &      315.5 &      266.0
& {\bf      340}  &     1103.5 &      492.5
 \\
  \hline
  \end{tabular}
 }

     {\footnotesize
      \tabcolsep 6pt
  \begin{tabular}{cccc|rrr|rrr|rrr}
  \hline
  \multicolumn{4}{c|}{{\small\bf hybrid} $\mNhyb^{-1}$        }  &  \multicolumn{3}{c|}{$p=1$} & \multicolumn{3}{c|}{$p=2$}&   \multicolumn{3}{c}{$p=3$}\\
  \hline
     $ \# \Th $   &
     $ N $ & $ \# \Thi $   &
     $ \# \ThH $ & 
     $ \iter $ & $ \MFl $ & $ \Mcomm $ & 
     $ \iter $ & $ \MFl $ & $ \Mcomm $ & 
     $ \iter $ & $ \MFl $ & $ \Mcomm $ \\ 
  \hline
    7549 &        7 &     1078 &        7
& {\bf      155}  &       36.2 &        9.9
& {\bf      233}  &      247.1 &       29.6
& {\bf      226}  &      853.9 &       47.9
 \\
   14937 &       14 &     1066 &       14
& {\bf      158}  &       39.4 &       27.0
& {\bf      208}  &      243.8 &       71.0
& {\bf      236}  &      812.8 &      134.2
 \\
   29819 &       29 &     1028 &       29
& {\bf      162}  &       39.1 &       70.4
& {\bf      215}  &      243.1 &      186.9
& {\bf      234}  &      861.0 &      339.0
 \\
  \hline
  \end{tabular}
 }

    \vspace{2mm}
    
    {
       domain decomposition {\bf respecting} component interfaces}
    
     {\footnotesize
      \tabcolsep 6pt
  \begin{tabular}{cccc|rrr|rrr|rrr}
  \hline
  \multicolumn{4}{c|}{ {\small\bf additive} $\mNadd^{-1}$     }  &  \multicolumn{3}{c|}{$p=1$} & \multicolumn{3}{c|}{$p=2$}&   \multicolumn{3}{c}{$p=3$}\\
  \hline
     $ \# \Th $   &
     $ N $ & $ \# \Thi $   &
     $ \# \ThH $ & 
     $ \iter $ & $ \MFl $ & $ \Mcomm $ & 
     $ \iter $ & $ \MFl $ & $ \Mcomm $ & 
     $ \iter $ & $ \MFl $ & $ \Mcomm $ \\ 
  \hline
    7549 &        7 &     1078 &        7
& {\bf      405}  &       80.5 &       25.7
& {\bf      644}  &      593.4 &       81.9
& {\bf      803}  &     2124.2 &      170.2
 \\
   14937 &       14 &     1066 &       14
& {\bf      405}  &       83.6 &       69.1
& {\bf      641}  &      613.5 &      218.7
& {\bf      738}  &     1950.2 &      419.7
 \\
   29819 &       29 &     1028 &       29
& {\bf      409}  &       92.8 &      177.7
& {\bf      609}  &      591.3 &      529.3
& {\bf      707}  &     2007.6 &     1024.2
 \\
  \hline
  \end{tabular}
 }

     {\footnotesize
      \tabcolsep 6pt
  \begin{tabular}{cccc|rrr|rrr|rrr}
  \hline
  \multicolumn{4}{c|}{{\small\bf hybrid} $\mNhyb^{-1}$        }  &  \multicolumn{3}{c|}{$p=1$} & \multicolumn{3}{c|}{$p=2$}&   \multicolumn{3}{c}{$p=3$}\\
  \hline
     $ \# \Th $   &
     $ N $ & $ \# \Thi $   &
     $ \# \ThH $ & 
     $ \iter $ & $ \MFl $ & $ \Mcomm $ & 
     $ \iter $ & $ \MFl $ & $ \Mcomm $ & 
     $ \iter $ & $ \MFl $ & $ \Mcomm $ \\ 
  \hline
    7549 &        7 &     1078 &        7
& {\bf      323}  &       65.0 &       20.5
& {\bf      516}  &      485.2 &       65.6
& {\bf      620}  &     1692.8 &      131.4
 \\
   14937 &       14 &     1066 &       14
& {\bf      317}  &       67.0 &       54.1
& {\bf      475}  &      473.0 &      162.1
& {\bf      519}  &     1454.8 &      295.2
 \\
   29819 &       29 &     1028 &       29
& {\bf      306}  &       72.5 &      133.0
& {\bf      425}  &      437.1 &      369.4
& {\bf      484}  &     1485.1 &      701.1
 \\
  \hline
  \end{tabular}
 }

  \end{center}
\end{table}

The observations for this more complicated case are, in principle, the same as those in
Section~\ref{sec:lapl}.
Neglecting the computational costs of the coarse solver, we have the weak scalability
of the algebraic solvers, cf. Tables~\ref{tab:alter2} and \ref{tab:alter5}.
The hybrid preconditioner $\mNhyb^{-1}$ saves about 25 -- 30\% of iterations compared to the additive preconditioner $\mNadd^{-1}$.
It is a potential benefit of the reduction of the computational costs,
if the computational costs of the coarse solver can be neglected.
On the other hand, if the coarse global problem is larger than the local ones,
the additive Schwarz method dominated the hybrid one, since the hybrid method cannot solve
the coarse problem in parallel.

Additionally, these observations are robust for the problem data. That is,
the increase of the ratio $\zeta=k_1/k_0 = \mu_0 /\mu_1$ (cf. \eqref{Kij} and \eqref{alter1})
by factor 100 leads to the increase of $\iter$ by 50 -- 75\% only. Hence, the convergence
of iterative solvers is significantly less sensitive to $\zeta$ compared to
the prediction resulting from the spectral bounds \eqref{RES1} -- \eqref{RES2}.

Finally, comparing the pairs of results corresponding
to the domain decomposition {\it respecting} and {\it non-respecting} component interfaces, we
found the latter one is more efficient only for the case $\zeta=\numf{10000}$ and
only when the coarse mesh has a sufficiently
large number of mesh elements $\#\ThH$.
We suppose that it is caused by the geometry of the problem. The domain decomposition
attempts to balance the number of degrees of freedom and the air component $\Om_a$ is relatively small,
so only a few coarse elements are defined in this component. Therefore, their shape is ``complicated''
and the corresponding coarse problem does not provide enough information.
This explanation is supported by the comparison in Section~\ref{sec:alterN}
where the air component is enriched due to mesh adaptation and the coarse element does not
suffer from previous problems. Then the domain decomposition respecting component interfaces
dominates.

\subsection{Nonlinear problem with $hp$-mesh adaptation}
\label{sec:alterN}

Finally, we employ the presented iterative method for the numerical solution
of a non-linear elliptic problem in combination with anisotropic $hp$-mesh
adaptation. Following \cite{Glowinski74,DolCon_JCAM23}, the problem geometry is the
same as in Figure~\ref{fig:alter_geom} and the magnetic potential $u$ fulfills 
\begin{align}
  \label{alt6}
  -\nabla\cdot \left(\nu (x, |\nabla u(x)|^2) \nabla u\right) &= f\qquad \mbox{ in }\Om
\end{align}
with \begin{align}
  \label{alt3}
  \nu(x, r) = 
  \begin{cases}
    \frac{1}{\mu_0} & \mbox{ for } x\in \Om_a, \\
    \frac{1}{\mu_0}\left(\alpha + (1-\alpha)\frac{r^4}{\beta+ r^4}\right) &
    \mbox{ for } x\in \Om_s\cup\Om_r. \\
  \end{cases}
\end{align}
The quantity
$\mu_0=1.256\times10^{-6}\, \mathrm{kg}\cdot\mathrm{m}\cdot\mathrm{A}^{-2}\cdot\mathrm{s}^{-2}$
denotes the permeability of
the vacuum and the material coefficients are $\alpha = 0.0003$, $\beta = 16000$ according
to \cite{Glowinski74}.
We consider the constant current density $f=5\times10^{4}\, \mathrm{A}\cdot\mathrm{m}^{-2}$ and
prescribe the mixed Dirichlet/Neumann boundary conditions as in Section~\ref{sec:alter}.
We note that these conditions differ from those in \cite{DolCon_JCAM23}.

We discretize \eqref{alt6} again using the SIPG 
method (cf. \cite{DolCon_JCAM23}), which leads to the similar form as
\eqref{Ah} but $\Ah$ is nonlinear in the first argument since $\K = \nu(x, |\nabla u(x)|^2)$
depends also on $u$. The arising nonlinear algebraic system is solved iteratively by
the Newton method where the Jacobian is evaluated by the differentiation 
of $\Ah$.

The iterative Newton method is stopped when the ratio between the algebraic error estimator
and the discretization error estimator is below $10^{-3}$, cf. \cite{hp-steady} for details.
At each Newton iteration, we solve a linear algebraic system by
the conjugate gradient method with additive and/or hybrid Schwarz preconditioners,
we employ the stopping criterion \eqref{alg:stop} with $\TOL=10^{-2}$.
The lower value of $\TOL$  leads to a small decrease in the number of Newton steps $\iterN$ but
a significant increase in the CG iterations $\iterL$.

\begin{figure}
  \includegraphics[height=0.30\textwidth]{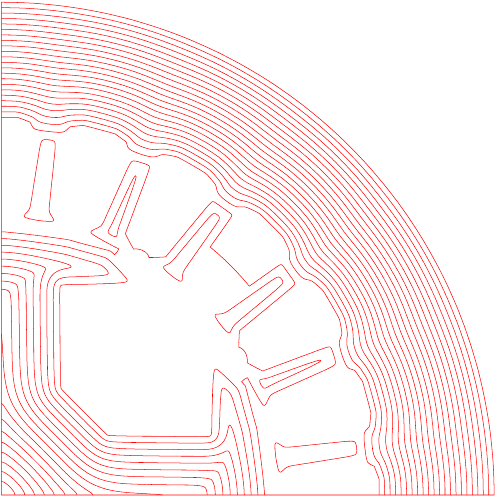}
  \hspace{0.01\textwidth}
  \includegraphics[height=0.30\textwidth]{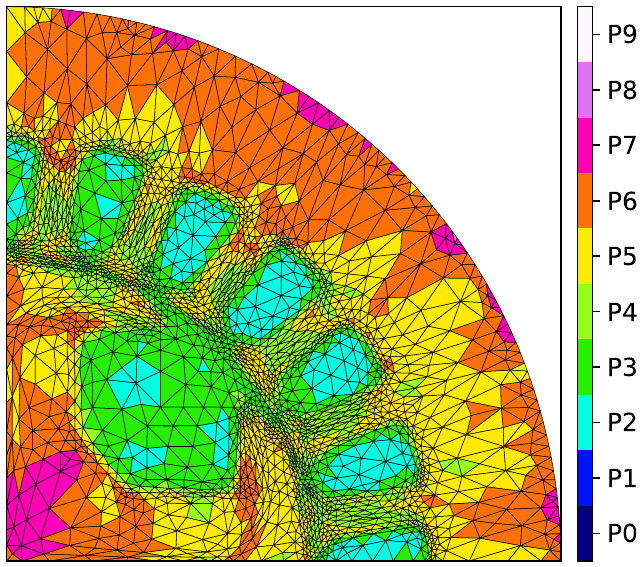}
  \hspace{0.01\textwidth}
  \includegraphics[height=0.30\textwidth]{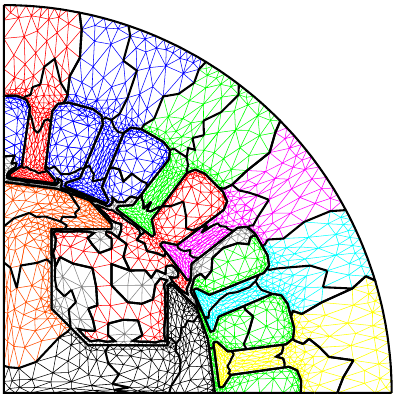}
  \caption{Nonlinear alternator \eqref{alt6}--\eqref{alt3}, isolines of the solution (left),
  the $hp$-mesh (center) and the domain decomposition (right) after 8~levels of mesh adaptation.}
  \label{fig:alterN} 
\end{figure}

Moreover, when the criterion for the Newton method is reached, we perform the re-meshing
using the anisotropic $hp$-mesh adaptation based on the interpolation error control, we refer
to \cite[Chapters~5-6]{AMA-book} for details. After the re-meshing (including the variation
of polynomial approximation degrees), a new domain decomposition of $\Om$ is employed and
the computational process is repeated.
We carried out eight levels of mesh refinement; at each adaptation level,
$\Om$ is divided into $N=12$ subdomains,
and the coarse mesh always has $\#\ThH=48$ elements. Figure~\ref{fig:alterN} shows the final
isolines of the solution, the final $hp$-mesh, and the final domain decomposition.
The subdomains $\Om_i$, $i=1,\dots,N$ are colored, and the coarse elements are plotted
by thick black lines.


Table~\ref{tab:alterN} presents
the comparison of the performance of the additive and hybrid Schwarz preconditioners on the
same sequence of adaptive meshes. The solution from the previous level is used as the initial
approximation of the solution on the next adaptive level.
For each level of mesh adaptation, we show the number of mesh elements $\#\Th$, the number
of degrees of freedom $\DoF=\dim\Shp$, and the number of
Newton iterations $\iterN$ and CG iterations $\iterL$.
In the last line of this table, we give the total number of $\MFl$ and $\Mcomm$.
We deduce a mild dominance of the hybrid preconditioner.
However, we are aware that a rigorous comparison requires a
deeper study and will be the subject of further research.

Furthermore, we observe
that the domain decomposition respecting the component interfaces strongly dominates.
As mentioned above, we suppose that this is
due to adapted meshes having more elements close to the interfaces, and then the
coarse elements have a better shape in narrow regions.

\begin{figure}
  \begin{center}
    \includegraphics[width=0.95\textwidth]{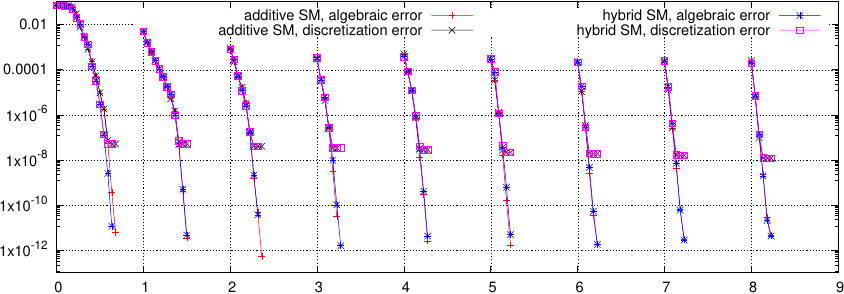}
  \end{center}
  \caption{Nonlinear alternator \eqref{alt6}--\eqref{alt3}, comparison of the convergence
    of the Newton method using CG solver with additive and hybrid Schwarz preconditioners
    (domain decomposition respecting the component interfaces).}
  \label{fig:alterNc}
\end{figure}

Finally, Figure~\ref{fig:alterNc} shows the convergence of the Newton method
using CG solver with additive and hybrid Schwarz preconditioners.
The horizontal axes correspond to the mesh adaptation level, and each node corresponds
to one Newton iteration. We plot the estimators of the algebraic and discretization errors
from \cite{hp-steady} for both preconditioners. We observe that space errors stagnate
in each mesh adaptation loop at the
level between $10^{-8}$ and $10^{-7}$ and the nonlinear solver stops when the
ratio between the algebraic and discretization estimators is $10^{-3}$.
\begin{table}
  \caption{Nonlinear alternator \eqref{alt6}--\eqref{alt3}, the convergence of the
    mesh adaptive algorithm.}
  \label{tab:alterN}
  \begin{center}

    {
       domain decomposition {\bf non-respecting} component interfaces}

     {
      \tabcolsep 4pt
  \begin{tabular}{ccc|rrrr|rrrr}
  \hline
  adaptive 
  & & &
  \multicolumn{4}{c|}{additive preconditioner $\mNadd^{-1}$} & \multicolumn{4}{c}{hybrid preconditioner $\mNhyb^{-1}$} \\
  \cline{4-11}
     level
 & $ \#\Th $ &  $ \DoF\ $ 
 & $ \iterN $ & $ \iterL\  $ &  $\ \MFl $ &$\ \Mcomm $   
 & $ \iterN $ & $ \iterL\  $  &  $\ \MFl $ &$\ \Mcomm $   
  \\
 \hline
 0 &  \numf{2154}  & \numf{12924}
 & \numf{18} & \numf{1817}  &  \numf{176.279999}  & \numf{84.1849976}
 & \numf{19} & \numf{1435} &  \numf{202.300003}  & \numf{66.4860001} \\ 
 1 &  \numf{1659}  & \numf{12153}
 & \numf{13} & \numf{1942}  &  \numf{274.190002}  & \numf{84.6090012}
 & \numf{14} & \numf{1664} &  \numf{308.350006}  & \numf{72.4970016} \\ 
 2 &  \numf{1688}  & \numf{14211}
 & \numf{9} & \numf{2172}  &  \numf{344.339996}  & \numf{110.650002}
 & \numf{9} & \numf{1802} &  \numf{366.290009}  & \numf{91.8050003} \\ 
 3 &  \numf{1838}  & \numf{17930}
 & \numf{8} & \numf{2653}  &  \numf{649.869995}  & \numf{170.529999}
 & \numf{8} & \numf{2299} &  \numf{666.280029}  & \numf{147.779999} \\ 
 4 &  \numf{2197}  & \numf{24301}
 & \numf{10} & \numf{3653}  &  \numf{1351.80005}  & \numf{318.239990}
 & \numf{8} & \numf{2379} &  \numf{1016.59998}  & \numf{207.250000} \\ 
 5 &  \numf{2609}  & \numf{32197}
 & \numf{9} & \numf{3490}  &  \numf{2177.89990}  & \numf{402.829987}
 & \numf{8} & \numf{2710} &  \numf{1840.19995}  & \numf{312.799988} \\ 
 6 &  \numf{3128}  & \numf{43462}
 & \numf{11} & \numf{4581}  &  \numf{4807.70020}  & \numf{713.760010}
 & \numf{8} & \numf{3067} &  \numf{3453.00000}  & \numf{477.869995} \\ 
 7 &  \numf{3845}  & \numf{60492}
 & \numf{11} & \numf{4619}  &  \numf{10066.0000}  & \numf{1001.70001}
 & \numf{9} & \numf{3588} &  \numf{8224.50000}  & \numf{778.099976} \\ 
 8 &  \numf{4589}  & \numf{82439}
 & \numf{13} & \numf{5558}  &  \numf{19735.0000}  & \numf{1642.59998}
 & \numf{10} & \numf{3979} &  \numf{14767.0000}  & \numf{1176.00000} \\ 
  \hline
    \multicolumn{3}{l}{total costs} & & & \numf{42233} & \numf{5205}
    & & & \numf{30844} & \numf{3330} \\

  \hline
    \end{tabular}
 }

    \vspace{5mm}
    
    {
       domain decomposition {\bf respecting} component interfaces}

     {
      \tabcolsep 4pt
  \begin{tabular}{ccc|rrrr|rrrr}
  \hline
  adaptive 
  & & &
  \multicolumn{4}{c|}{additive preconditioner $\mNadd^{-1}$} & \multicolumn{4}{c}{hybrid preconditioner $\mNhyb^{-1}$} \\
  \cline{4-11}
     level
 & $ \#\Th $ &  $ \DoF\ $ 
 & $ \iterN $ & $ \iterL\  $ &  $\ \MFl $ &$\ \Mcomm $   
 & $ \iterN $ & $ \iterL\  $  &  $\ \MFl $ &$\ \Mcomm $   
  \\
 \hline
 0 &  \numf{2154}  & \numf{12924}
 & \numf{16} & \numf{1362}  &  \numf{129.240005}  & \numf{63.1040001}
 & \numf{15} & \numf{1038} &  \numf{140.520004}  & \numf{48.0929985} \\ 
 1 &  \numf{1670}  & \numf{12307}
 & \numf{12} & \numf{1138}  &  \numf{138.020004}  & \numf{48.4510002}
 & \numf{12} & \numf{828} &  \numf{138.809998}  & \numf{35.2519989} \\ 
 2 &  \numf{1710}  & \numf{14417}
 & \numf{9} & \numf{1822}  &  \numf{328.489990}  & \numf{90.8720016}
 & \numf{8} & \numf{1055} &  \numf{243.649994}  & \numf{52.6180000} \\ 
 3 &  \numf{1867}  & \numf{17747}
 & \numf{6} & \numf{657}  &  \numf{198.550003}  & \numf{41.7999992}
 & \numf{7} & \numf{708} &  \numf{253.910004}  & \numf{45.0449982} \\ 
 4 &  \numf{2214}  & \numf{23986}
 & \numf{7} & \numf{1029}  &  \numf{461.369995}  & \numf{88.4830017}
 & \numf{7} & \numf{794} &  \numf{427.760010}  & \numf{68.2750015} \\ 
 5 &  \numf{2653}  & \numf{32569}
 & \numf{6} & \numf{1317}  &  \numf{864.489990}  & \numf{153.770004}
 & \numf{6} & \numf{920} &  \numf{727.739990}  & \numf{107.419998} \\ 
 6 &  \numf{3209}  & \numf{45266}
 & \numf{6} & \numf{1756}  &  \numf{1707.19995}  & \numf{284.959991}
 & \numf{6} & \numf{1298} &  \numf{1448.30005}  & \numf{210.639999} \\ 
 7 &  \numf{3915}  & \numf{63482}
 & \numf{6} & \numf{1558}  &  \numf{3553.30005}  & \numf{354.570007}
 & \numf{6} & \numf{1308} &  \numf{3297.69995}  & \numf{297.679993} \\ 
 8 &  \numf{4798}  & \numf{88242}
 & \numf{6} & \numf{1744}  &  \numf{7536.29980}  & \numf{551.700012}
 & \numf{6} & \numf{1373} &  \numf{6680.10010}  & \numf{434.339996} \\ 
  \hline
    \multicolumn{3}{l}{total costs} & & & \numf{14916} & \numf{1677}
    & & & \numf{13358} & \numf{1299} \\

  \hline
    \end{tabular}
 }

  \end{center}
\end{table}

\section{Conclusion}
\label{sec:concl}

We presented the two-level hybrid Schwarz preconditioner for the $hp$-discontinuous
Galerkin discretization of elliptic problems as an alternative to
the well-established additive Schwarz approach. Whereas numerical analysis shows
only a slightly better spectral bound, numerical experiments
demonstrate about 25 -- 30\% savings in the number of iterations of the linear solver used.
The savings in the number of communication operations are similar,
but the savings in the number of floating point operations per core are smaller since the coarse
problem cannot be solved in parallel to the local ones.
Therefore, the weak scalability is shown for both preconditioners
provided that the computational costs of the coarse solver are negligible.
Furthermore, robustness with respect to the data problem, anisotropy of the meshes
and varying polynomial degree is demonstrated.
Finally, we showed that domain decomposition taking into account material interfaces
dominates for fine meshes when the coarse elements have reasonable shape.

The future research will orient to more challenging problems (nonlinear, nonsymmetric,
systems of equations). Additionally, the optimal choice of subdomains and coarse grid
deserve a deeper insight in order to balance the number of floating-point operations
and communication operations.

\section*{Acknowledgments}
We are thankful to prof. Martin Gander (University of Geneva)
for a fruitful discussion during his stay at the Charles University in Prague and
our colleagues Michal Outrata and Petr Tich\'y inspiring suggestions and ideas.

 \section*{Conflict of interest}

 The authors declare that they have no conflict of interest.

 \section*{Data Availability Statement}

 Data sets were not generated or analyzed during the current study.


\end{document}